\documentclass
{amsart}
\usepackage{amsmath, amsthm, amssymb,mathrsfs,amsfonts,}
\usepackage{mathtools}
\usepackage{relsize}
\usepackage{enumitem}
\usepackage[utf8]{inputenc}
\usepackage{pgffor}

\usepackage
{hyperref}

\usepackage{mathrsfs}

\usepackage{cases} 
\usepackage{tikz}

\usepackage[normalem]{ulem}

\newcommand{\stkout}[1]{\ifmmode\text{\sout{\ensuremath{#1}}}\else\sout{#1}\fi}

\newcommand{\Be}{\begin{equation}}
\newcommand{\Ee}{\end{equation}}
\newcommand{\Bea}{\begin{eqnarray}}
\newcommand{\Eea}{\end{eqnarray}}
\newcommand{\Bel}{\begin{align}}
\newcommand{\Eel}{\end{align}}
\newcommand{\Beas}{\begin{eqnarray*}}
	\newcommand{\Eeas}{\end{eqnarray*}}
\newcommand{\Benu}{\begin{enumerate}}
	\newcommand{\Eenu}{\end{enumerate}}
\newcommand{\Bi}{\begin{itemize}}
	\newcommand{\Ei}{\end{itemize}}

\renewcommand{\a}{\bar a}
\numberwithin{equation}{section}
\newcommand{\dist} {\text{dist\! }}

\theoremstyle{plain}
\newtheorem{thm}{Theorem}[section]
\newtheorem{cor}[thm]{Corollary}

\newtheorem{lem}[thm]{Lemma}
\newtheorem{prop}[thm]{Proposition}

\newtheorem{defn}[thm]{Definition}

\definecolor{ao}{rgb}{0, 0.5, 0}

\newcommand{\sz}{s_{\!\,\mathsmaller{0}}}

\newcommand{\fa}{\mathfrak a}

\newcommand{\otheta}{{\bar{\theta}}}

\newcommand{\mt}{\mathfrak m_t}
 
\newcommand{\fN}{\mathfrak N}

\newcommand{\ssum}[1]{\sum{}_{{}_{{\mathlarger #1}}}}

\DeclareMathOperator{\Vol}{Vol}
\DeclareMathOperator{\sspan}{span}
\DeclareMathOperator{\supp}{supp}

\newcommand{\ndeg}[2]{\fN(#1, #2)}
\newcommand{\vcon}[2]{\fV(#1, #2)}

\newcommand{\fV}{\mathfrak V}
\newcommand{\sza}{s_{\mathsmaller{\!\star}}}
\newcommand{\delc}{\delta_{\mathsmaller{\!\star}}}
\newcommand{\szz}{s_{\!\,\mathsmaller{0}}} 
\newcommand{\bfa}{\bar{\mathfrak a}}

\begin{document}

\title[Sharp smoothing of averages  over curves]{
Sharp smoothing properties of  \\  averages  over curves 
}

\author[H. Ko]{Hyerim Ko}
\author[S. Lee]{Sanghyuk Lee}
\author[S. Oh]{Sewook  Oh}

\address{Department of Mathematical Sciences and RIM, Seoul National University, Seoul 08826, Republic of Korea}
\email{kohr@snu.ac.kr}
\email{shklee@snu.ac.kr}
\email{dhtpdnr0220@snu.ac.kr}


\begin{abstract}    
We prove sharp smoothing  properties  of  the averaging operator defined by convolution with a  measure  on a smooth nondegenerate 
curve $\gamma$ in  $\mathbb R^d$, $d\ge 3$.  Despite the simple geometric structure of such curves,  
the sharp smoothing estimates  have remained largely unknown except for those  in low dimensions. 
Devising a novel inductive strategy, we obtain  the optimal  $L^p$ Sobolev regularity  estimates,  which settle the conjecture raised by Beltran-Guo-Hickman-Seeger.  Besides,  we show the sharp local smoothing estimates for every $d$. 
As a result,  we establish, for the first time, nontrivial $L^p$ boundedness of the maximal average over dilations of $\gamma$ for $d\ge 4$. 
\end{abstract}

\maketitle 

\section{Introduction}
The regularity property  of  integral transforms defined by averages over submanifolds
is a fundamental  subject  in harmonic analysis, which  has been extensively studied since the 1970s.  
There is an immense body of literature devoted to the subject (see, for example,  \cite{SW, PhS, Stein93, CNSW}  and references therein). However, numerous problems remain wide open.     
The regularity  property is typically addressed in the frameworks  of  $L^p$ improving, $L^p$ Sobolev regularity,  and local smoothing estimates, to which 
 $L^p$ boundedness of  the maximal average is also closely related.  In this paper, we study the smoothing estimates for the averaging operator given by convolution with a measure supported on a curve. 
\thispagestyle{empty}

Let $I=[-1,1]$ and $\gamma $ be a smooth curve from $I$ to $\mathbb R^d$.   
We define a measure $\mathfrak m_t$ supported on $t\gamma$   by  
\[ \langle  \mathfrak m_t , f \rangle  = \int  f(t\gamma(s)) \psi(s) ds,   \] 
where $\psi\in \mathrm C_c^\infty((-1,1))$.  
We are concerned  with $d\ge 3$ since all the problems we address in the current paper  are  well understood  when  $d=2$. 
We consider the averaging  operator 
\[
\mathcal A_t f(x)= f\ast \mt(x)
\]
and study  the above-mentioned  regularity problems on $\mathcal A_t $ under the assumption that $\gamma$ is nondegenerate, 
that is to say, 
 \begin{align}\label{nonv}
\det( \gamma'(s), \dots, \gamma^{(d)}(s) )\neq0,
\quad  s\in I.
\end{align}

The $L^p$ improving property of  $ \mathcal A_t$ for a fixed $t\neq 0$ now has a complete characterization,  see \cite{Christ2, Stovall} (also, see \cite{TW} for  generalizations to variable coefficient settings). However,  $L^p$ Sobolev and local smoothing estimates  for  $\mathcal A_t$  turned out to be more involved  and  are far  less well understood.   Recently,  there has been progress  in low dimensions  $d= 3, 4$ (\cite{PS, KLO, BGHS, BGHS2}), but it does not seem feasible to extend the approaches in   recent works to higher dimensions. We discuss this  matter  in  detail near the end of the introduction.
 By devising  an inductive strategy,  we prove the optimal  $L^p$ Sobolev regularity and sharp local  smoothing estimates  in any dimension $d\ge 3$.  
  As a result, we also obtain $L^p$ boundedness of 
 the associated maximal function which was unknown for $d\ge 4$.

\subsubsection*{$L^p$ Sobolev regularity.} Let  $2\le p\le \infty$. We set  
$
\mathcal A f=\mathcal A_1 f
$
and consider the $L^p$ Sobolev regularity  estimate
\begin{align}\label{A1}
\| \mathcal A f\|_{L_\alpha^p(\mathbb R^d)}
\le C \|f\|_{L^p(\mathbb R^d)}.
\end{align}
When $d=2$,  the estimate holds if and only if $\alpha\le 1/p$ (e.g., see \cite{Christ}). 
In higher dimensions, however,  the problem of obtaining \eqref{A1} with the sharp smoothing order $\alpha$ becomes highly nontrivial except for 
the $L^2\to L^2_{1/d}$ estimate which is an easy consequence of the decay  of Fourier transform of  ${\mathfrak m}_t$:
 \Be \label{f-decay}  |\widehat {\mathfrak m}_t(\xi)|\le C(1+|t\xi|)^{-1/d}\,.\Ee  
 
 It was conjectured by Beltran, Guo, Hickman, and Seeger 
\cite[Conjecture 1]{BGHS2} that \eqref{A1}  holds  for  $\alpha\le 1/p$ if $2d-2<p<\infty$. 
When $d=3$, the conjecture  was verified by 
 the conditional result of Pramanik and Seeger \cite{PS}  and the decoupling inequality due to  Bourgain and Demeter  \cite{BD} (see \cite{OSS, TV} for earlier results).  
 The case $d=4$ was recently obtained  by Beltran et al \cite{BGHS2}.  Our first result  proves the conjecture for every  $d\ge 5$. 
 
 \begin{thm}
\label{thm:improving}
Let $d\ge3$. Suppose $\gamma$ is a smooth  nondegenerate curve.  Then, \eqref{A1}  holds for $\alpha\le 1/p$ if $p> 2(d-1)$.
\end{thm}

 Interpolation with the $L^2\to L^2_{1/d}$ estimate gives \eqref{A1} for  $\alpha<(p+2)/(2dp)$ when $2< p \le 2(d-1)$. 
 It is also known that   \eqref{A1}  fails  if $\alpha> \alpha(p):=\min (1/p,  (p+2)/(2dp)) $ (see  \cite[Proposition 1.2]{BGHS2}). 
 Thus,  only  the  estimate \eqref{A1} with $\alpha=\alpha(p) $ remains open  for $2< p \le 2(d-1)$.  Those endpoint estimates seem to be 
 a subtle problem.  The argument in this paper provides simpler alternative proofs  of the previous  results for $d=3,4$.  Theorem  \ref{thm:improving}  remains valid as long as $\gamma\in  \mathrm C^{2d}(I)$ (see Theorem \ref{sob-induct}).   However, we do not try to optimize the regularity assumption.
  
 The result in Theorem \ref{thm:improving} can be easily generalized  to curves of different types.   We say a smooth curve $\gamma$ from $I$ to $\mathbb R^d$ is  of finite  type  if  there is an $ \ell$ such that  $\sspan \{ \gamma^{(1)}(s), \dots, \gamma^{( \ell)}(s)\}=\mathbb R^d$ for each $s\in I$.  The  type  at $s$ is defined to be the smallest of such $ \ell$ and the maximal type   is the supremum over $s \in I$ of the type at $s$.  (See, e.g., \cite{PS, HL}.)  Using Theorem \ref{thm:improving} and a rescaling argument  (\cite{PS, HL})  one can obtain   the following, which  proves the Conjecture 2 in  \cite{BGHS2}.

\begin{cor}\label{type}
Let $d\ge3$, $ \ell> d$ and $2 \le p <\infty$. Suppose $\gamma$ is a curve of maximal type $ \ell$.
Then \eqref{A1} holds for $\alpha \le \min \big(\alpha(p) , \,1/ \ell\big)$  if  $p\neq \ell$ when $ \ell\ge 2d-2$, and if  $p\in  [2,  2\ell/(2d-\ell))\cup ( 2d-2, \infty)$  when $d< \ell< 2d-2$. 
\end{cor}

By interpolation \eqref{A1} holds for $\alpha < \min \big(\alpha(p) , \,1/ \ell\big)$   if $p= \ell$ when $ \ell\ge 2d-2$, and  if ${2 \ell}/{(2d- \ell)} \le  p \le 2d-2$ when $d< \ell< 2d-2$.
These estimates are sharp. Since a finite type curve contains a nondegenerate subcurve  and the   $L^2\to L^2_{1/ \ell}$ estimate  is optimal, 
\eqref{A1} fails if $\alpha > \min \big(\alpha(p) , \,1/ \ell\big)$.  When $ \ell\ge 2d-2$, Corollary \ref{type} resolves  the problem of 
the Sobolev regularity estimate \eqref{A1}.  In fact, the failure of $L^ \ell\to L^ \ell_{1/ \ell}$ bound was shown in  \cite{BGHS2} using Christ's example \cite{Christ}.  
By \cite[Theorem 1.1]{SeW},  Corollary  \ref{type} also gives $H^1(\mathbb R^d)\to  L^{1,\infty}(\mathbb R^d)$ bound on the lacunary maximal function
$f\to \sup_{k\in \mathbb Z} |f\ast \mathfrak m_{2^k}|$ 
  whenever $\gamma$ is of finite type.

\subsubsection*{Sharp local smoothing}  
We now consider the  estimate  
\begin{align}
\label{ls}
\| 
\chi(t) \mathcal A_t f\|_{L_\alpha^p(\mathbb R^{d+1})}
\le C \|f\|_{L^p(\mathbb R^d)},
\end{align}
where  
$\chi$ is a smooth function supported in  $(1/2,4)$. Compared with the   $L^p$ Sobolev estimate  \eqref{A1},  
the additional integration in $t$ is expected to yield extra smoothing. Such a phenomenon is called \emph{local smoothing}, which  
has been studied for the dispersive equations to a great extent (e.g., see \cite{Sjolin, CS}). However, the local smoothing for the averaging operators  exhibits  considerably different nature. 

In particular,  there is no local smoothing when $p=2$. Besides, a bump function example shows $\alpha\le 1/d$.  As we shall see later,  
the estimate \eqref{ls} fails unless  $\alpha\le 2/p$ (Proposition \ref{prop:n}).  So,  it seems to be plausible to conjecture that 
\eqref{ls} holds for $\alpha<\min(2/p, 1/d)$ if $2<p<\infty$.  
For $d=2$, the conjecture follows  by  the recent  result on Sogge's local smoothing conjecture  for the wave operator (\cite{Sogge, Wolff, LV, BD}), which is due to Guth,  Wang, and Zhang 
\cite{GWZ}.  
 When $d=3$, some local smoothing estimates 
were utilized  by Pramanik and Seeger \cite{PS} and  Beltran et al. \cite{BGHS}  to prove $L^p$ maximal bound.

Nevertheless,  for $d\ge 3$, no local smoothing estimate  up to  the sharp order $2/p$ 
has been known previously.  

\begin{thm}\label{LS}
Let $d\ge 3$. Suppose $\gamma$ is a smooth nondegenerate curve. Then,  if $p\ge 4d-2$,  
\eqref{ls} holds true for  $\alpha<2/p$. 
\end{thm}

 Theorem \ref{LS} remains valid as far as  $\gamma\in \mathrm C^{3d+1}(I)$ (see Theorem \ref{lclb} below).  
 
 \subsubsection*{Maximal estimate}  The local smoothing estimate \eqref{ls} has been of particular interest in connection to   $L^p$ boundedness of 
 the maximal operator
\begin{align*}
M f(x)=\sup_{0<t} |\mathcal A_t f(x)| 
\end{align*}
(\cite{MSS, SS, PS, BGHS}) 
and problems in geometric measure theory (see, e.g.,  \cite{Wolff} and Corollary \ref{packing} below).  
If  the estimate \eqref{ls} holds for some $\alpha>1/p$, $L^p$  boundedness of  $M$  follows  by a standard argument relying on  the Sobolev embedding (\cite{PS}).

The study of the maximal functions generated by dilations of submanifolds goes back to Stein's spherical maximal theorem \cite{Stein}  
(see, also, \cite[Ch.10]{Stein93} and \cite{IKM}).  
The circular maximal theorem was later proved by 
Bourgain  \cite{B2} (also, see \cite{Sogge,  MSS,  Schlag97,   SS,  L}).  Afterwards, a natural question 
was whether the maximal operator $M$ under consideration in the current paper  is bounded on  $L^p$ for some  $p\neq \infty$ when $d\ge 3$.  In view of Stein's  interpolation  argument  based on $L^2$ estimate (\cite{Stein}),  proving  $L^p$ boundedness of $M$   becomes more challenging  as $d$ increases since 
the decay of  the Fourier transform of $\mt$  gets weaker (see \eqref{f-decay}).     Though the question  was raised as early as  in the late 1980s, it  remained open for any $d\ge 3$ until recently.  
In  $\mathbb R^3$, 
the first positive result was obtained  by Pramanik and Seeger  \cite{PS}  and the range of $p$  was further extended to $p>4$ thanks to 
the decoupling inequality for the cone \cite{BD}.   
Very recently, the authors \cite{KLO} proved $L^p$ boundedness of $M$ on the optimal range, i.e., $M$ is bounded on $L^p$ if and only if $p>3$. 
The same result was independently obtained by Beltran et al. \cite{BGHS}. 
 
However, no nontrival $L^p$ bound on $M$ has been known in higher dimensions.   
The following establishes existence of 
such  maximal  bound for every $d\ge 4$.

\begin{thm}\label{max}
Let $d\ge4$. Suppose $\gamma$ is  a smooth nondegenerate curve.
Then, for  $p>2(d-1)$  we have 
\begin{align}\label{Mf}
\| Mf \|_{L^p(\mathbb R^d)} \le C \|f\|_{L^p(\mathbb R^d)}.
\end{align}
\end{thm}

The result is a  consequence of Theorem \ref{LS}.  
Since the estimate  \eqref{ls} holds for $p=2$ and $\alpha=1/d$,  interpolation 
gives  \eqref{ls} for  some $\alpha>1/p$ when $2d-2<p<\infty$.  So,  the maximal estimate \eqref{Mf} follows, as mentioned before,  by  a standard argument.  
A natural conjecture is that  $M$ is bounded on $L^p$   if and only if $p>d$.  $M$ 
can not be bounded on $L^p$ if $p\le d$, as can be seen  by a simple adaptation of the argument in \cite[Proposition 4.4]{KLO}.   Theorem \ref{max}  also extends to the  finite type curves by a rescaling argument (\cite{HL, PS}).  The following result  is sharp when $\ell \ge 2(d-1)$.

\begin{cor} 
Let $d\ge4$ and $\ell> d$. Suppose $\gamma$ is a curve of maximal type $\ell$.
Then \eqref{Mf} holds if $p>\max(\ell,2(d-1))$.
\end{cor}

\subsubsection*{Packing of curves in $\mathbb R^d$} 
The sharp local smoothing estimate \eqref{ls} in Theorem \ref{LS}  has  interesting measure theoretic consequences  concerning unions of 
curves generated  by  translation and dilation of a nondegenerate curve. The following generalizes  Wolff's result  \cite[Corollary 3]{Wolff}, where unions of  circles in $\mathbb R^2$ were considered (see also \cite{Marstrand, Mitsis, Wolff97} for earlier results).

\begin{cor}\label{packing}
Let $\gamma$ be a smooth nondegenerate curve in $\mathbb R^d$, $d\ge 3$, and let $E \subset \mathbb R^{d+1}$ be a set of Hausdorff dimension greater than $d-1$. 
Suppose  $F$ is a set in $\mathbb R^d$ such that  $(x+t\gamma(I))\cap F$ has positive $1$-dimensional outer measure for  all $(x,t) \in E$. Then $F$ has positive outer measure. 
\end{cor}

Corollary \ref{packing} follows by  Theorem \ref{LS} and  the argument in \cite{Wolff}.  
The result does not hold in general without the nondegeneracy assumption on $\gamma$ as one can easily  see  considering  a curve  contained in a lower dimensional affine  space. 
The same result continues to be valid for the finite type curve.  Consequently, Corollary \ref{packing} implies the following. 
 
\medskip
 
\noindent {\bf Corollary  \ref{packing}$'$\!.} \emph{Let $\gamma$ be a smooth finite type curve in $\mathbb R^d$, $d\ge 3$, and let $E$ and $F$  be  compact subsets in $\mathbb R^d$. 
Suppose  $E$ has Hausdorff dimension greater than $d-1$ and   for each  $x \in E$ there is $t(x)> 0$ such  that $x+t(x)\gamma(I)\subset F$. Then, $F$ has positive measure.}

\subsubsection*{Our approach}  
To prove $L^p$ ($p\neq  2$) smoothing properties of $\mathcal A_t$, we need more than the decay of  $\widehat {\mathfrak m}_t$, i.e., \eqref{f-decay}. 
 When $d=2$, we have rather a precise  asymptotic expansion  of $\widehat {\mathfrak m}_t$, which makes it possible to relate $\mathcal A_t$ to 
 other forms of operators.  In fact, one can use the estimate for the wave operator (e.g., \cite{SS, TV, L})  to obtain local smoothing estimate.  However,  in higher dimensions $d\ge 3$,   to compute $\widehat {\mathfrak m}_t$ explicitly is not  a simple matter. 
     Even worse, this becomes much  more complicated  as $d$ increases since one has to take into account  the derivatives $\gamma^{(k)}(s)\cdot\xi$, $k=2,\dots, d$. 
 The common  approach  in  \cite{PS, BGHS, BGHS2}  to get around this difficulty  was  to use detailed decompositions (of various scales)  on the Fourier side
away from the conic sets where  $\widehat {\mathfrak m}_t$ decays slowly.    
The consequent decompositions were then  combined with the decoupling or square function estimate \cite{OSS, PS1, PS, PS2, BGHS, BGHS2}.  
 However, this type of  approach based on fine scale decomposition  becomes exceedingly difficult  to manage as the dimension $d$  gets larger  and, consequently,   
does not seem to be tractable in higher dimensions.

To overcome the difficulty,  we develop a new strategy  which allows us to dispense with such sophisticated decompositions. 
 Before closing the introduction we briefly discuss the key ingredients of our approach. 

\noindent $\bullet$ 
The main novelty of the paper lies in an  induction argument which we build on the local nondegeneracy assumption: 
\Be
 \label{sumN}
 \tag*{$\mathfrak N(L, B)$}  \sum_{\ell=1}^L |\langle \gamma^{(\ell)}(s),  \xi\rangle| \ge B^{-1}|\xi|  
 \Ee  
 for a constant  $B\ge 1$.
  To prove our results, we consider the operator $\mathcal A_t[\gamma, a]$ (see \eqref{Aa} below for its definition).
 Clearly, $\ndeg{d}{B'}$ holds for a constant $B'>0$  if  $\gamma$ satisfies  \eqref{nonv}. 
 However, instead of considering the  case  $L=d$ alone,    we prove the estimate for all $L=2,\dots, d$ under the assumption  that \ref{sumN} holds on the support of $a$. See Theorem \ref{lclb} and \ref{sob-induct}.  
A trivial (yet, important)  observation is that  $\ndeg{L-1}B$ implies $\ndeg LB, $
so we may think  of  $\mathcal A_t[\gamma, a]$ as being more degenerate as $L$ gets larger.    
 Thanks to this  hierarchical structure, we may use an inductive strategy  along the number $L$. 
See Proposition \ref{newprop} and  \ref{newpropsob} below.   

\noindent $\bullet$ We extend   the rescaling  \cite{HL, KLO}  and iteration \cite{PS} arguments. Roughly speaking, we combine the first with the induction assumption in Proposition \ref{newprop} (or  \ref{newpropsob}) to handle the less degenerate parts, and use the latter to deal with the remaining part. 
In order to generalize the arguments,  we introduce  a class of symbols which are naturally adjusted to a small subcurve  (Definition \ref{symbol}). 
We also use  the decoupling inequalities for the nondegenerate curves obtained  by  Beltran et al. \cite{BGHS2} (Corollary \ref{rcdecoupling}). 
Their  inequalities were deduced from those due to Bourgain, Demeter, and Guth \cite{BDG}. 
Instead of applying the inequalities directly, we use modified forms which are adjusted to  the sharp smoothing orders of the specific estimates (see \eqref{dgain0} and \eqref{dgain}). 
This  makes it possible  to obtain the sharp estimates  on an extended range.  

 \smallskip
\noindent\emph{Organization of the paper.} 
We first prove Theorem \ref{LS} whose proof is more involved than that of  Theorem \ref{thm:improving}. 
In Section \ref{sec:ls}, we reduce the proof of Theorem \ref{LS}  to 
that of Proposition \ref{microlocal}, 
which we prove while assuming Proposition \ref{iterative}.  The proof of  Proposition \ref{iterative} is given in  Section  \ref{sec:mainthm}. 
We prove Theorem \ref{thm:improving} in Section \ref{sec:sobolev}.

\section{Smoothing estimates with localized frequency} \label{sec:ls}
In this section, we consider an extension of  Theorem \ref{LS}  via microlocalization  (see Theorem \ref{lclb} below) 
which we can prove inductively.   We then  reduce the matter to proving Proposition \ref{microlocal}, which we show by applying Proposition \ref{iterative}.
We also  obtain some preparatory results.

 Let $1\le L \le d$ be a positive  integer and $B\ge1$  be a large number. 
For quantitative control of estimates we consider the following two conditions:   
\begin{align}
\label{curveB}
 &\qquad \max_{0\le j\le 3d+1}  |\gamma^{(j)}(s)|\,\le\,  B,    && s\in I, 
\\
\label{lindepN}
\tag*{$\mathfrak V(L,B)$}
&\Vol \big( \gamma^{(1)}(s), \dots,\gamma^{(L)}(s) \big)\ge 1/B,    &&s\in I,
\end{align} 
where $ \Vol (v_1,\dots,v_L)$ denotes  the $L$-dimensional  volume of the parallelepiped generated by 
 $v_1,\dots, v_L$.  By finite decomposition, rescaling, and a change of variables,  the constant $B$ can be taken to be close to $1$ (see  Section \ref{Rescaling}).

 \medskip 
\noindent{\bf Notation.} For nonnegative quantities $A$ and $ D$,  we denote  $A\lesssim D$ if there exists an independent positive  constant $C$ such that $A\le CD$, 
but  the constant   $C$ may differ at each occurrence  depending  on the context, and $A \lesssim_{B}\! D$ means the inequality holds with an implicit constant depending on $B$.  
Throughout the paper, the constant $C$ mostly depends on $B$. However, we do not make it explicit every time since it is clear in the context. 
By  $A=O(D)$  we  denote $|A|\lesssim D$.

\begin{defn}
For $k\ge 0$, we  denote $
\mathbb A_k=\{ \xi \in \mathbb R^d: 2^{k-1} \le |\xi| \le 2^{k+1}\}.$  
We say $a\in  \mathrm C^{d+L+2}(\mathbb R^{d+2})$
  is a  symbol of type $(k, L, B)$ relative to  $\gamma$ if 
$\supp a $ $\subset I \times [2^{-1},4] \times \mathbb A_k$,  
\ref{sumN}
holds  for $\gamma$ whenever $(s,t,\xi)\in \supp a$ for some $t$,  and     
\begin{align*}
|\partial_s^{j}\partial^{l}_t\partial_\xi^\alpha a(s,t,\xi)|\le  B |\xi|^{-|\alpha|} \end{align*}
for $(j,l,\alpha)\in \mathcal I_L:=\{(j,l,\alpha):  0\le j\le1,\, 0\le l\le 2L,\, |\alpha|\le d+L+2\}$. 
 \end{defn}
We define  an integral operator by
\begin{align}
\label{Aa} 
\mathcal A_t[\gamma, a]f(x)=(2\pi)^{-d}\iint_{\mathbb R}  e^{i(x-t\gamma(s))\cdot\xi}\,a(s,t,\xi)ds \,\widehat f (\xi)\,d\xi.
\end{align}
Note $\mathcal A_t f=\mathcal A_t[\gamma, \psi]f$.  
Theorem \ref{LS} is a consequence of  the following.

\begin{thm}\label{lclb} 
Let  $\gamma\in \mathrm C^{3d+1}(I)$ satisfy \eqref{curveB} and \ref{lindepN} for some   $B\ge 1$. Suppose $a$ is a symbol of type $(k, L, B)$ relative to $\gamma$. 
 Then,   if  $p\ge 4L-2$,  for any $\epsilon>0$  there is a constant $C_\epsilon=C_\epsilon(B)$ such that 
\begin{align}\label{smoothing}
\|
\mathcal A_t[\gamma, a]f\|_{L^p(\mathbb R^{d+1})}\le
C_\epsilon 2^{(-\frac{2}p+\epsilon) k}\|f\|_{L^p(\mathbb R^d)}. 
\end{align}
\end{thm}

Theorem  \ref{lclb}  is trivial when $L=1$. Indeed, \eqref{smoothing} follows from the estimate $|\mathcal A_t[\gamma, a]f(x)|\lesssim_B   \int_I K\ast |f|(x-t\gamma(s)) \,ds$ where $K(x)= 2^{(d-1)k}(1+ |2^{k}x|)^{-d-3}.$ To show this, note $|\gamma'(s)\cdot\xi|\sim 2^k$ if $(s,t,\xi)\in \supp a$ for some $t$. By
integration by parts in $s$, $\mathcal A_t[\gamma, a]=t^{-1}\mathcal A_t[\gamma, \tilde a]$ where $\tilde a=i(\gamma'(s)\cdot\xi\, \partial_s a -\gamma''(s)\cdot\xi\, a)/(\gamma'(s)\cdot\xi)^2$. 
Since $|\partial_\xi^{\alpha}  \tilde a|\lesssim |\xi|^{-|\alpha|-1}$ for $|\alpha|\le d+3$,  
routine integration by parts in $\xi$ gives the estimate (e.g., see \emph{Proof of Lemma \ref{kernel}}). 
When $L=2$,  Theorem  \ref{lclb} is already known  by the result  in \cite[Theorem 4.1]{PS} and  the decoupling inequality in \cite{BD}.

Once we have Theorem  \ref{lclb}, the proof of Theorem  \ref{LS} is straightforward.  By Littlewood-Paley decomposition  it is sufficient to show  
 \eqref{smoothing}  for $p\ge 4d-2$ with $a_k(s,t,\xi)=\psi(s) \chi(t) \beta(2^{-k}|\xi|)$, where $\beta \in \mathrm C_c^\infty((1/2,2))$. This can be made rigorous using $\iint  e^{-i t(\tau+\gamma(s)\cdot\xi)}   \psi(s) \chi(t) dsdt$ $=O((1+|\tau|)^{-N})$ for any $N$ if $|\tau|\ge(1+\max_{s\in \supp \psi} |\gamma(s)|)|\xi|$.   
 Since  $\gamma$  satisfies \eqref{nonv}, $a_k$ is of type $(k, d, B)$ relative to  $\gamma$ for a large $B$. Therefore, Theorem \ref{LS} follows from Theorem \ref{lclb}. 
  
Theorem \ref{lclb} is  immediate  from the next proposition, which places Theorem \ref{lclb} in 
an inductive framework. 

\begin{prop}\label{newprop}
Let $2 \le N  \le d$. Suppose Theorem \ref{lclb} holds for $L=N-1$.
Then, Theorem \ref{lclb} holds true with $L=N$.
\end{prop}

To prove Proposition \ref{newprop}, 
from this section to Section \ref{sec:mainthm} we fix $N\in [2, d]$, $\gamma$ satisfying 
$\vcon N B$, and  a symbol 
$ a$  of type $(k, N, B)$ relative to $\gamma$.

 One of the main ideas  is that  by a suitable decomposition of the symbol 
 we can separate  from  $\mathcal A_t[\gamma,a]$ the less degenerate part   which corresponds to $L=N-1$. 
 To this part we apply the   assumption combined with a rescaling argument.  To do this,  we   introduce 
 a class of symbols  which are adjusted to short subcurves of $\gamma$. 

\subsection{Symbols associated to subcurves}  
We begin with some notations.  Let  $N\ge2$, and  let $\delta$ and $B'$  denote the numbers such that 
\begin{align*}
2^{-k /N}\le \delta\le 2^{-7dN}B^{-6N},    \qquad  B \le  B' \le B^C
\end{align*}
for a large constant $C\ge3d+1$.  We note that  $\vcon{N-1}{B'}$ holds  for some $B'$.  
In fact,  $\vcon{N-1}{B^2}$  follows  by  \eqref{curveB} and $\vcon NB$. 

For $s\in I$,  we define a linear map $ \widetilde{\mathcal L}^\delta_s:\mathbb R^{d}\mapsto \mathbb R^{d}$  
as follows: 
\begin{equation}
\begin{aligned}\label{overL}
(\widetilde{\mathcal L}_{s}^\delta)^\intercal \gamma^{(j)}(s)
&=\delta^{N-j}\gamma^{(j)}(s), \qquad  && j=1,\dots,N-1,
\\
(\widetilde{\mathcal L}_{s}^\delta)^\intercal v &=v, \qquad \qquad   \quad~ && v\in  \big(\mathrm V_s^{\gamma,N-1}\big)^{\perp},
\end{aligned}
\end{equation} 
where $\mathrm V_s^{\gamma,\ell} = \sspan \big\{ \gamma^{(j)}(s) : j=1,\dots,\ell\big\}$. $\widetilde {\mathcal L}_{s}^\delta$ is well-defined
since $\vcon{N-1}{B^2}$ holds for $\gamma$. The linear map $\widetilde {\mathcal L}_{s}^\delta$ naturally appears when we rescale a subcurve of length about $\delta$ (see  the  proofs of Lemma \ref{kernel} and \ref{lem:res}). 
We denote 
\Be \label{L}
\mathcal L_{s}^\delta(\tau, \xi)=
\big(\delta^N \tau - \gamma(s)\cdot   \widetilde{\mathcal L}_s^\delta \xi,\,\, \widetilde{\mathcal L}_s^\delta \xi\big),  \qquad (\tau, \xi)\in \mathbb R\times \mathbb R^d.
\Ee   

We set $G(s)=(1,\gamma(s))$ and define 
\begin{align*}
\Lambda_k(s,\delta, B')= \bigcap_{0\le j\le  N-1}
\big\{(\tau,\xi) \in \mathbb R \times \mathbb A_k & : 
|\langle G^{(j)}(s),(\tau,\xi)\rangle| \le  B'2^{k+5} \delta^{N-j}
\big
\}. 
\end{align*}

\begin{defn} 
\label{symbol} 
Let  $(\sz,\delta)\in (-1,1)\times(0, 1)$ such that $I(\sz,\delta):=[\sz-\delta,  \sz +\delta]\subset I$.  Then, by $\mathfrak A_k (\sz,\delta)$ 
we denote the set of  $\mathfrak a \in \mathrm C^{d+N+2}(\mathbb R^{d+3})$
such that
\begin{align}\label{symsupp}
 & \supp  \mathfrak a  \subset    I(\sz,\delta)\times [2^{-1}, 2^2]\times \Lambda_k(\sz,\delta, B),
  \\
\label{symineq2}
&\big|\partial^{j}_s\partial^{l}_t\partial^{\alpha}_{\tau, \xi}
\mathfrak a\big( s,t,  \mathcal L_{\sz}^\delta(\tau,\xi) \big) \big|
\le B\delta^{-j}|(\tau,\xi)|^{-|\alpha|}, \qquad (j,l, \alpha)\in \mathcal I_N. 
\end{align}
\end{defn}

We define $\supp_{\xi} \mathfrak a=\bigcup_{s,t,\tau} \supp \mathfrak a(s,t,\tau,\cdot)$ and $\supp_{s, \xi} \mathfrak a=\bigcup_{t,\tau} \supp \mathfrak a(\cdot,t,\tau,\cdot)$,  
and  $\supp_{s} \mathfrak a$ and $\supp_{\tau, \xi} \mathfrak a$ are defined  likewise.  We note  \emph{a  statement $S(s, \xi)$, depending on  $(s, \xi)$, 
 holds on   $\supp_{s,\xi} \mathfrak a$  if and only  if $S(s,\xi)$ holds whenever $(s,t, \tau,\xi) \in \supp \mathfrak a$ for some $t$, $\tau$}.

\medskip

Denote $\mathrm V_s^{G,\ell} = \sspan \{(1,0),  G'(s), \dots,  G^{(\ell)}(s) \}.$ We take a close look at the map ${\mathcal L}_s^\delta$. 
By \eqref{overL} and \eqref{L} we have 
\Be
\label{hoho} 
\begin{cases}
\begin{aligned}
(\mathcal L_s^\delta)^\intercal G(s) &= \delta^{N}(1,0),
\\[2pt]
(\mathcal L_s^\delta)^\intercal  G^{(j)}(s) &= \delta^{N-j}G^{(j)}(s), \qquad &&j=1,\dots, N-1,  
\\
({\mathcal L}_s^\delta)^\intercal v &=v,    \qquad&& v\in  (\mathrm V_s^{G,N-1} )^\perp.  
\end{aligned}
\end{cases}
\Ee 
The first identity is clear since $(\mathcal L_s^\delta)^\intercal(\tau, \xi)=(\delta^N \tau, (\widetilde{\mathcal L}_s^\delta)^\intercal \xi-\tau (\widetilde{\mathcal L}_s^\delta)^\intercal \gamma(s))$. 
The second and the third follow from \eqref{overL} since  $G^{(j)}\in \{0\} \times \mathbb R^d$, $1 \le j \le N-1$, $\big(
\mathrm V_s^{G,N-1} 
\big)^\perp\subset \{0\} \times \mathbb R^d$, and $(\mathcal L_s^\delta)^\intercal(0,\xi)=(0,(\widetilde{\mathcal L}_s^\delta)^\intercal \xi)$.   
Furthermore,  there is  a constant $C=C(B)$, independent of $s$ and $\delta$, such that 
\Be
\label{Lnorm}
| \mathcal L_{s}^{\delta} (\tau, \xi)|  \le C|(\tau,\xi)|.
\Ee
Note that \eqref{Lnorm} is equivalent to  $ |( \mathcal L_{s}^{\delta})^\intercal (\tau, \xi)|  \le C|(\tau,\xi)|$. The inequality  is clear  from 
\eqref{overL}  because  $\vcon{N-1}{B^2}$  holds 
and all the eigenvalues of $(\widetilde{\mathcal L}_{s}^\delta)^\intercal$ are contained in the interval $(0,1]$.

\begin{lem}\label{2ksupport}  
Let $ \mathcal L_{s}^\delta(\tau,\xi)\in\Lambda_k(s,\delta, B')$ and $\vcon{N-1}{B'}$ holds for $\gamma$. Then, there exists  a constant $C=C(B')$ such that 
\begin{align}\label{C2}
C^{-1} |(\tau,\xi)|\le 2^k \le C|\xi|.
\end{align}
\end{lem}
\begin{proof}
Since $ \mathcal L_{s}^\delta(\tau,\xi)\in\Lambda_k(s,\delta, B')$,  by \eqref{L}  we have $2^{k-1}\le|\widetilde{\mathcal L}_{s}^\delta\xi|\le 2^{k+1}$.  
So, the second inequality in \eqref{C2} is clear from  \eqref{Lnorm}  if we take  $\tau=0$.

To show the first inequality, from \eqref{hoho} we have 
$| \langle (1,0),(\tau,\xi) \rangle | \le B'2^{k+5}$ and
 $|\langle G^{(j)}(s),(\tau, \xi) \rangle| \le B' 2^{k+5}$, $1\le j\le N-1$, because $\mathcal L_{s}^\delta(\tau,\xi)\in \Lambda_k(s,\delta, B').$ Also, if  $v \in ( \mathrm V_s^{G,N-1}  )^\perp$  
and  $|v|=1$,  by \eqref{hoho}
 we see $ |\langle v,(\tau,\xi) \rangle| =|\langle v, \mathcal L_{s}^\delta(\tau,\xi) \rangle|\le 2^{k+1}$.  Since $\vcon{N-1}{B'}$ holds and $\mathrm V_s^{G,N-1}\oplus ( \mathrm V_s^{G,N-1} )^\perp=\mathbb R^{d+1}$,  
we get $|(\tau,\xi)| \le C2^k$ for some $C=C(B')$. 
\end{proof}

The following shows the matrices $ \mathcal L_{s}^{\delta}$, $\mathcal L_{\sz}^{\delta}$ are close to each other if  so are $s, \sz$.  

\begin{lem}\label{similarmx}
Let  $s, \sz \in (-1,1)$ and  $\gamma$ satisfy $\vcon{N-1}{B'}$. 
If $|s-\sz |\le \delta $, then 
there exists a constant $C=C(B')\ge1$ such that
\begin{align}\label{LLt}
C^{-1}|(\tau,\xi)| \le | (\mathcal L_{\sz}^{\delta})^{-1} \mathcal L_{s}^{\delta} (\tau, \xi)|  \le C|(\tau,\xi)|.
\end{align} 
\end{lem}

\begin{proof}
It suffices to prove that \eqref{LLt} holds if $|s-\sz| \le  c\delta$ for a   constant $c>0$, 
independent of $s$ and $\sz$. Applying this  finitely many times, we can remove the additional assumption.  
Moreover,  it is enough to show 
\Be 
\label{L-L}
\|(\mathcal L_{s}^{\delta})^{\intercal} (\mathcal L_{\sz}^\delta)^{-\intercal}- \mathrm I\| \lesssim_{B'} c
\Ee
when $|s-\sz| \le  c\delta$.  Here,  $\| \cdot\|$ denotes a matrix norm. 
Taking $c>0$ sufficiently small, we get \eqref{LLt}. 

By \eqref{hoho},   $(\mathcal L_{s}^{\delta})^{\intercal} (\mathcal L_{\sz}^\delta)^{-\intercal}G^{(j)}(\sz)
= \big( \mathcal L_{s}^\delta\big)^{\intercal} \delta^{-(N-j)}G^{(j)}(\sz)$ for $j=1,\dots,N-1$. Let $\sz=s+c'\delta$, $|c'|\le c$.   
Expanding $G^{(j)}$ in  Taylor series  at $s$, by \eqref{curveB} we have
\begin{align*}
(\mathcal L_{s}^{\delta})^{\intercal} (\mathcal L_{\sz}^\delta)^{-\intercal}G^{(j)}(\sz)
&=\big( \mathcal L_{s}^\delta\big)^{\intercal}
\Big( \sum_{\ell=j}^{N-1} \delta^{-(N-j)}G^{(\ell)}(s) \frac{(
c'\delta)^{\ell-j}}{(\ell-j)!}  + O\big(c^{N-j}B'\big)\Big)
\end{align*}
for $j=1,\dots,N-1$. By \eqref{hoho} and   the mean value theorem, we get 
\[
(\mathcal L_{s}^{\delta})^{\intercal} (\mathcal L_{\sz}^\delta)^{-\intercal}G^{(j)}(\sz) =G^{(j)}(\sz)+O(cB'), \quad j=1,\dots,N-1. 
\]
From \eqref{hoho} we also have  $(\mathcal L_{s}^{\delta})^{\intercal} (\mathcal L_{\sz}^\delta)^{-\intercal}(1,0)=\delta^{-N} (\mathcal L_{s}^{\delta})^{\intercal} G(\sz)$. A similar argument also shows  $(\mathcal L_{s}^{\delta})^{\intercal} (\mathcal L_{\sz}^\delta)^{-\intercal}(1,0)=(1,0)+O(cB')$. 

Let  $\{ v_N,\dots, v_d\}$ denote an orthonormal basis of   $(  \mathrm V_{\sz}^{G,N-1} )^\perp$. 
By $\vcon{N-1}{B'}$ and \eqref{curveB} it follows that $|\gamma^{(j)}(\sz)|\ge (B')^{-1-N}$, $j=1, \dots, N-1$.
  Since $|\gamma^{(j)}(s)-\gamma^{(j)}(\sz)|\le\! cB'\delta$, there is an orthonormal basis $\{v_N(s),\dots,v_d(s)\}$ of
$(  \mathrm V_s^{G,N-1})^\perp$ such that $|v_j(s)-v_j| \lesssim_{B'}\! c\delta$, $j=N,\dots,d$. 
So,  we have 
$|(\mathcal L_{s}^\delta)^\intercal v_j -v_j|\lesssim_{B'}\! c\delta$ by \eqref{Lnorm}. 
Since $(\mathcal L_{\sz}^\delta)^{-\intercal} v_j = v_j$,  it follows that $|(\mathcal L_{s}^{\delta})^{\intercal} (\mathcal L_{\sz}^\delta)^{-\intercal} v_j -v_j|\lesssim_{B'}\! c\delta$, $j=N,\dots,d$.  

We denote by $\mathrm M$  the matrix  $[(1,0),G'(\sz), \dots, G^{(N-1)}(\sz), v_N, \dots, v_d]$. Then,  combining all together,   we have $\|(\mathcal L_{s}^{\delta})^{\intercal} (\mathcal L_{\sz}^\delta)^{-\intercal} \mathrm M-\mathrm M\|\lesssim_{B'}\! c$.  
Note that $\vcon{N-1}{B'}$ gives  $|\mathrm M^{-1} v|  \lesssim_{B'}\!  |v|$ for $v\in\mathbb R^{d+1}$.  Therefore,  we  obtain 
\eqref{L-L}. 
\end{proof}

 For a continuous function $\mathfrak a$ supported in $I\times [1/2, 4]\times \mathbb R\times \mathbb A_k$, we set 
\begin{align}\label{Ma}
m[\mathfrak a](\tau,\xi)&=\iint e^{-it'(\tau+ \gamma(s)\cdot \xi)}\mathfrak a(s,t',\tau,\xi)dsdt', 
\\[2pt]
\label{Ta}
\mathcal T[\mathfrak a]f(x,t)&=(2\pi)^{-d-1}  \iint  e^{i(x\cdot\xi+t\tau)}m[\mathfrak a](\tau,\xi) \widehat f (\xi)\,d\xi d\tau. 
\end{align}

\begin{lem}\label{kernel} 
Suppose $\mathfrak a\in C^{d+3}(\mathbb R^{d+3})$  satisfies  \eqref{symsupp} 
and \eqref{symineq2} for $j=l=0$ and $|\alpha| \le d+3$.  Then, there is a constant $C=C(B)$ such that 
\begin{align}
\label{ker-est} \| \mathcal T[\mathfrak a]f\|_{L^\infty(\mathbb R^{d+1})} &\le C\delta \|f\|_{L^\infty(\mathbb R^d)},
\\
\label{kernelrmk}
\| (1-\tilde\chi)\mathcal T[\mathfrak a]f\|_{L^p(\mathbb R^{d+1})} &\le C2^{-k}\delta^{1-N} \|f\|_{L^p(\mathbb R^d)}, \quad p> 1, 
\end{align}
where   $\tilde\chi\in \mathrm C_c^\infty((2^{-2},2^3))$ such that $\tilde\chi=1$ on $[3^{-1},6]$. 
\end{lem}

\begin{proof}  
 We first note   
\Be\label{conv}
 \mathcal T[\mathfrak a]f (x,t) = \int K[\mathfrak a](s,t, \cdot)\ast f(x) \,ds,
 \Ee
where 
\Be 
\label{eq:kernel}
K[\mathfrak a](s,t,x)=(2\pi)^{-d-1}  \iiint e^{i(t-t',x-t'\gamma(s))\cdot (\tau, \xi)}\mathfrak a(s,t',\tau,\xi) \, d\xi d\tau dt'. 
\Ee

Since  $\supp_s \mathfrak a \subset I(\sz, \delta)$, to prove   \eqref{ker-est} we need only  to show
\Be
 \label{ker1}
   \|K[\mathfrak a](s,\cdot)\|_{L^{\infty}_{t}L_x^1}\le C,  \qquad s\in I(\sz, \delta)
\Ee
for some $C=C(B)>0$.  To this end,  changing variables $(\tau,\xi) \rightarrow  2^k \mathcal L_{s}^{\delta}(\tau,\xi)$ in the right hand side of \eqref{eq:kernel} and noting  $|\!\det \mathcal L_{s}^{\delta}|= \delta^N |\!\det  \widetilde{\mathcal L}_{s}^\delta |$ $=\delta^{{N(N+1)}/2}$,   we get   
\begin{align*}
K[\mathfrak a](s,t,x)=C_\ast  \iiint e^{i 2^k(t-t',x-t\gamma(s))\cdot (\delta^N \tau,\, \widetilde{\mathcal L}_{s}^\delta \xi)} \mathfrak a(s,t',2^k\mathcal L_{s}^{\delta}(\tau,\xi)) \,  d\xi d\tau dt',
\end{align*}
where $C_\ast =(2\pi)^{-d-1} \delta^{{N(N+1)}/2} 2^{k(d+1)}$.  Since $\mathfrak a$ satisfies \eqref{symsupp},  
by   \eqref{LLt} and Lemma \ref{2ksupport} we have  $\supp\, \mathfrak a(s,t,2^k\mathcal L_{s}^\delta\cdot)\subset \{(\tau, \xi):  |(\tau,\xi)| \lesssim_B\! 1\}$.
Besides,   by  \eqref{symineq2} and \eqref{LLt} it follows that   $|\partial_{\tau,\xi}^{\alpha} \big(\mathfrak a(s,t, 2^k\mathcal L_{s}^\delta(\tau,\xi))\big)| \lesssim_B \!1 $ for $|\alpha|\le d+3$.  Thus, 
repeated integration by parts in $\tau,\xi$  yields
\[|K[\mathfrak a](s,t,x)| \lesssim C_\ast
\int_{1/2}^{4} \Big(1+2^k \big|\big (\delta^N(t-t'), (\widetilde{\mathcal L}_{s}^\delta)^\intercal(x-t\gamma(s))\big)\big| \Big)^{-d-3}\,dt',\]
by which  we obtain  \eqref{ker1} as desired.  

It is easy to show  \eqref{kernelrmk}.  The above estimate for $K[\mathfrak a]$ gives
\[\|(1-\tilde\chi)K[\mathfrak a](s,t,\cdot)\|_{L^1_x} \lesssim \delta^{-N} 2^{-k}|t-1|^{-1}|1-\tilde\chi(t)|.\]
Since $\supp_s \mathfrak a \subset I(\sz, \delta)$,   \eqref{kernelrmk} for $p>1$ follows   
by \eqref{conv}, Minkowski's and  Young's convolution inequalities. 
\end{proof}

\subsection{Rescaling} 
\label{Rescaling}
Let $\mathfrak a \in \mathfrak A_k(\sz,\delta)$. Suppose   that   
\begin{align}\label{lowD}
\sum_{j=1}^{N-1} \delta^{j} |\langle \gamma^{(j)}(s), \xi \rangle |\ge
{2^k \delta^N}/{B'} 
\end{align}
holds on $\supp_{s,\xi} \fa$ for some $B'>0$.  Then, via decomposition and  rescaling, we can bound the $L^p$ norm of  $\mathcal T[\mathfrak a]f$ by  
those of the operators given by symbols of type $(j, N-1, \tilde B)$ relative to a curve  for some $\tilde B$ and $j$ (see Lemma \ref{lem:res} below).  

To do so, 
we  define  a rescaled curve  $\gamma_{\sz}^\delta:I\to \mathbb R^d$ by
\begin{equation}
\label{tildegamma}
\gamma_{\sz}^\delta(s)
= \delta^{-N} (\widetilde{\mathcal L}_{\sz}^\delta)^\intercal 
\big(\gamma(\delta s+\sz)-\gamma(\sz)\big).
\end{equation}
As $\delta\to 0$, the curves $\gamma_{\sz}^\delta$ get close to a nondegenerate curve in $N$ dimensional vector space, 
so the curves behave in a  uniform way. In particular, \eqref{curveB} and $\vcon N B$ hold for some $B$ for $\gamma_{\sz}^\delta$ if $\delta <\delta'$ for a constant  $\delta'=\delta'(B)$ small enough. 

 Note $(\gamma_{\sz}^\delta)^{(j)}(s)=\delta^{j -N} (\widetilde{\mathcal L}_{\sz}^\delta)^\intercal\gamma^{(j)}(\delta s+\sz)$, $1\le j\le  N-1$, 
 and  $|( \gamma_{\sz}^\delta)^{(j)}(s)|\lesssim B\delta$,  $N+1\le j\le 3d+1$. Thus, Taylor series expansion and \eqref{overL} give 
\[(\gamma_{\sz}^\delta)^{(j)}(s)= \sum_{k=0}^{N-j-1}\frac{\gamma^{(j+k)}(\sz)}{k!}s^k+  \frac{(\widetilde{\mathcal L}_{\sz}^\delta)^\intercal \gamma^{(N)}(\sz)}{(N-j)!} s^{N-j} +O(B\delta)  \]  
for $j=1, \dots, N-1$.  By  \eqref{tildegamma}, we have  $(\gamma_{\sz}^\delta)^{(N)}(s)=(\widetilde{\mathcal L}_{\sz}^\delta)^\intercal \gamma^{(N)}(\sz)+O(\delta)$.
We write $\gamma^{(N)}(\sz)=v+v'$ where $v\in  \mathrm V_s^{\gamma,N-1}$ and $v'\in (\mathrm V_s^{\gamma,N-1})^{\perp}.$
Then, $(\widetilde{\mathcal L}_{\sz}^\delta)^\intercal \gamma^{(N)}(\sz) $ $= (\widetilde{\mathcal L}_{\sz}^\delta)^\intercal v+ v'$. Since 
$|(\widetilde{\mathcal L}_{\sz}^\delta)^\intercal v|\lesssim_B \delta$ and $|v'|\le B$, $|(\widetilde{\mathcal L}_{\sz}^\delta)^\intercal \gamma^{(N)}(\sz)|\le B+ C\delta$ for some 
$C=C(B)$.  Thus,  $\gamma=\gamma_{\sz}^\delta$ satisfies \eqref{curveB}  with $B$ replaced by $3B$ 
 if $\delta <\delta'$.

  An elementary argument (elimination) shows 
\begin{align*}
\Vol \big( (\gamma_{\sz}^\delta)^{(1)}(s),\dots,(\gamma_{\sz}^\delta)^{(N)}(s)  \big)
=\Vol \big( \gamma^{(1)}(\sz),\dots, \gamma^{(N)}(\sz) \big)+O(\delta)
\end{align*}
since $(\widetilde{\mathcal L}_{\sz}^\delta)^\intercal \gamma^{(N)}(\sz)= (\widetilde{\mathcal L}_{\sz}^\delta)^\intercal v+ v'$ and $\gamma^{(N)}(\sz)=v+v'$. 
Taking $\delta'$ small enough, from $\vcon NB$ for $\gamma$ we see $\vcon{N}{3B}$ hold for $ \gamma=\gamma_{\sz}^{\delta}$ if $0<\delta<\delta'$.

The next lemma (cf. \cite[Lemma 2.9]{KLO}) plays a crucial role in what follows.

\begin{lem}\label{lem:res}
Let $2\le N\le d$, $\mathfrak a \in \mathfrak A_k(\sz, \delta)$,  and $j_\ast=\log (2^k \delta^N)$. 
Suppose \eqref{lowD} holds on $\supp_{s,\xi} \fa$. 
Then, there exist  constants $C$, $\tilde B\ge 1$, and $\delta'>0$ depending on $B$, and  symbols 
$a_{1}, \dots, a_{l_\ast}$ of type $(j, N-1, \tilde B)$  relative to $\gamma_{\sz}^\delta$,  
such that 
\begin{align*}
\big\| \tilde\chi\, \mathcal T [\fa] f \big\|_{L^p(\mathbb R^{d+1})}
\le C \delta \sum_{1\le l\le C}  
\big\| \mathcal A_t 
[\gamma_{\sz}^{\delta}, \raisebox{-.2ex}{$a_l$} \,]
\tilde f_l \big\|_{L^p(\mathbb R^{d+1})}, 
\end{align*}
 $ \|\tilde f_l\|_p= \|f\|_p$,  
 and  $j\in [j_\ast-C,  j_\ast +C]$    
as long as  $0<\delta<\delta'$.
\end{lem}

\begin{proof} We set ${\mathfrak a}_{\delta, \sz}(s,t, \tau, \xi)=\mathfrak a\big(\delta s+\sz,t,\tau,\xi\big).$ Combining \eqref{Ma} and \eqref{Ta},  we write 
$\mathcal T[\mathfrak a] f$ as an integral (e.g.,  see \eqref{conv} and \eqref{eq:kernel}). Then, the change of  variables $s \rightarrow \delta s +\sz$ and $(\tau,\xi) \rightarrow (\tau-\gamma(\sz)\cdot \xi,\, \xi)$ gives  
\[
\mathcal T[\mathfrak a] f(x,t)=(2\pi)^{-d-1}  \, \delta
\iint e^{i \langle x-t\gamma(\sz),\, \xi \rangle} \mathcal J (s,t,\xi) \widehat f(\xi)\,ds d\xi,
\]
where
\[
\mathcal J(s,t,\xi)
=\!\!\iint e^{it \tau}e^{-it'(\, \tau+(\gamma(\delta s+\sz)-\gamma(\sz) )\cdot \xi\,)}
\, {\mathfrak a}_{\delta, \sz}(s,t',\tau-\gamma(\sz)\cdot \xi,\,\xi\big) \,dt' d\tau.
\] 

Let $\tilde f$ be given  by  $ \mathcal F(\tilde f)=|\det   \delta^{-N} \widetilde{\mathcal L}_{\sz}^\delta|^{1-1/p} \widehat f( \delta^{-N} \widetilde{\mathcal L}_{\sz}^\delta \cdot\,)$  where $\mathcal F( \tilde f\,)$ denotes the Fourier transform of $ \tilde f.$ Then, $ \|\tilde f\|_p=\|f\|_p$. Changing   variables  $\xi \rightarrow  \delta^{-N} \widetilde{\mathcal L}_{\sz}^\delta\xi$ gives 
\begin{align*}
\mathcal T [\mathfrak a] f (x,t)
= C_d
\iint e^{i \langle x-t\gamma(\sz),  \delta^{-N} \widetilde{\mathcal L}_{\sz}^\delta \xi \rangle}  
\mathcal J (s,t, \delta^{-N} \widetilde{\mathcal L}_{\sz}^\delta \xi) \mathcal F (\tilde f\,)(\xi)
\,dsd\xi, 
\end{align*} 
where $C_d=(2\pi)^{-d-1}   \,
\delta |\! \det  \delta^{-N} \widetilde{\mathcal L}_{\sz}^\delta |^{1/p}$.  This leads us to  set 
\begin{align}\label{tildeA}
\tilde a(s,t,\xi)= \frac1{2\pi} 
\iint e^{-it'(\tau+\gamma_{\sz}^\delta(s)\cdot\xi)}
\tilde\chi(t) {\mathfrak a}_{\delta, \sz}\big(s, t+t',\delta^{-N}\mathcal L_{\sz}^\delta(\tau,\xi) \big) \,dt' d\tau.
\end{align}

 It is easy to check $\tilde a\in C^{d+N+2}(\mathbb R^{d+2})$, since so is  $\mathfrak a$ and $\gamma\in \mathrm C^{3d+1}$. 
By \eqref{tildegamma} and \eqref{L}, we note $\tilde\chi(t)\mathcal J(s,t,  \delta^{-N} \widetilde{\mathcal L}_{\sz}^\delta \xi)={2\pi}  e^{-it  \gamma_{\sz}^\delta(s)\cdot\xi}\, \tilde a(s,t,\xi)$. Therefore,  
\begin{align*} 
\tilde\chi(t)\mathcal T [\fa]f(x,t)
=\delta |\det  \delta^{-N} \widetilde{\mathcal L}_{\sz}^\delta|^{\frac 1p}\, \mathcal A_t 
[ {\gamma_{\sz}^\delta}, \raisebox{-.3ex}{$\tilde a$} \,]
\tilde f\, \big( \delta^{-N} (\widetilde{\mathcal L}_{\sz}^\delta)^{\intercal} (x-t\gamma(\sz))\big),
\end{align*}
and a change of  variables  gives 
\Be\label{res:TA1}
\big\| \tilde\chi\, \mathcal T [\fa] f \big\|_{L^p(\mathbb R^{d+1})}
=  \delta
\big\| \mathcal A_t 
[\gamma_{\sz}^{\delta}, \raisebox{-.2ex}{$\tilde a$} \,]
\tilde f \big\|_{L^p(\mathbb R^{d+1})}.
\Ee

We shall obtain symbols of type $(j, N-1,\tilde B)$ from $\tilde a$ via decomposition and rescaling. To this end, we first  note  
\Be\label{fourierF}   \supp_\xi \tilde a \subset \big\{ \xi  \in \mathbb R^d: C^{-1}\delta^N 2^k \le |\xi | \le C \delta^N 2^{k} \big\} \Ee
for a constant $C=C(B)\ge 1$. 
This follows by Lemma \ref{2ksupport}  since 
 there exists $\tau$ such that  $\delta^{-N}\mathcal L_{\sz}^\delta(\tau,\xi) \in \Lambda_k(\sz, \delta, B)$ if $\xi \in \supp_\xi \tilde a$.  
 We  claim 
 \Be
\label{diff:a}
 |\partial_s^{j} \partial^{l}_t\partial_\xi^\alpha \tilde a(s,t,\xi)  | \lesssim_B\!  |\xi|^{-|\alpha|}, \qquad (j,l,\alpha)\in \mathcal I_{N-1}. 
\Ee

To show \eqref{diff:a}, let us set
\begin{align*}
\mathfrak b(s,t,t',\tau,\xi)
=\tilde\chi(t)  {\mathfrak a}_{\delta, \sz}(s, t+ t',\delta^{-N}\mathcal L_{\sz}^\delta(\tau,\xi)  \big).
\end{align*}
Note $0\le j\le 1$. Taking derivatives on both sides of \eqref{tildeA}, we have
\[
\partial_s^{j}\partial^{l}_t \partial_\xi^\alpha \tilde a(s,t,\xi) =\mathcal I[\mathfrak b_1]:
=\,\frac1{2\pi} \iint e^{-it'(\tau+\gamma_{\sz}^\delta(s)\cdot\xi)}
\mathfrak b_1(s,t,t',\tau,\xi) \,dt' d\tau, 
\]
where 
\begin{align*}
\mathfrak b_1=\!\!
\sum_{\substack {u_1+u_2=j, \\\alpha_1+\alpha_2+\alpha_3=\alpha}} \!\!
C_{\alpha, u}\big(t' \gamma_{\sz}^{\delta\,\,\prime}\cdot \xi \big)^{u_1-|\alpha_1|}
(t'\gamma_{\sz}^{\delta\,\,\prime})^{\alpha_1}
(t'\gamma_{\sz}^\delta)^{\alpha_2}
\, \partial_{s}^{u_2}\partial^{l}_t \partial_{\xi}^{\alpha_3} \mathfrak b,
\end{align*}
with $0 \le u_1 \le 1$,  $0 \le |\alpha_1| \le u_1$, and constants $C_{\alpha, u}$  satisfying  $|C_{\alpha, u}|=1$. 
Integration by parts  $u_1+|\alpha_2|$ times in $\tau$   gives  $\partial_s^{j}\partial^{l}_t\partial_\xi^\alpha \tilde a=\mathcal I[\mathfrak b_2]$, 
where
\begin{align*}
\mathfrak b_2 
= 
 \!\!\!\sum_{\substack {u_1+u_2=j, \\\alpha_1+\alpha_2+\alpha_3=\alpha}} \!\!\!
C_{\alpha, u}' \big(\gamma_{\sz}^{\delta\,\,\prime}\cdot \xi \big)^{u_1-|\alpha_1|}
(\gamma_{\sz}^{\delta\,\,\prime})^{\alpha_1}
(\gamma_{\sz}^\delta)^{\alpha_2} \,	\partial_\tau^{u_1+|\alpha_2|}  \partial_s^{u_2} \partial^{l}_t \partial_{\xi}^{\alpha_3} \mathfrak b
\end{align*}
with  constants  $C_{\alpha, u}'$ satisfying  $|C_{\alpha, u}'|=1$. We decompose  $\mathcal I[\mathfrak b_2]= \mathcal I[\chi_E \mathfrak b_2]+  \mathcal I[\chi_{E^c}\mathfrak b_2]$ where 
$E=\{(\tau, \xi): |\tau+ \gamma_{\sz}^\delta(s)\cdot \xi|\le 1\}$.  Then, integrating by parts in $t'$ for $ \mathcal I[\chi_{E^c}\mathfrak b_2]$,  we obtain 
\[
|\mathcal I[\mathfrak b_2]| \lesssim   \iint  \chi_E |\mathfrak b_2| \,+   \frac{\chi_{E^c} |\partial_{t'}^2 \mathfrak b_2|}{|\tau+ \gamma_{\sz}^\delta(s)\cdot \xi|^{2}}  
\,dt'd\tau.
\]
Since $\mathfrak a \in \mathfrak A_k(\sz,\delta)$, 
$
|\partial_s^{j'} \partial_t^{l'} \partial_{\tau,\xi}^{\alpha'} \mathfrak b| \lesssim_B\! 
|\xi|^{-|\alpha'|}
$
for $(j', l', \alpha')\in \mathcal I_N$.  It is also clear that $|\gamma_{\sz}^{\delta\,\,\prime}(s)|\lesssim 1$ if  $\delta< \delta' $. Thus,  $|\mathfrak b_2|=O(|\xi|^{-|\alpha|})$ and $|\partial_{t'}^2\mathfrak b_2|=O(|\xi|^{-|\alpha|})$ if $l\le 2(N-1) $.  Since $\partial_s^{j}\partial^{l}_t\partial_\xi^\alpha \tilde a=\mathcal I[\mathfrak b_2]$, we get  \eqref{diff:a}.

Now, we decompose $\tilde a$. Let $\tilde \chi_1, \tilde\chi_2$,  and $\tilde\chi_3\in C_c^\infty(\mathbb R)$  such that  $\tilde \chi_1 +\tilde\chi_2+\tilde\chi_3 =1 $ on $\supp \tilde\chi$ and $\supp\tilde \chi_\ell \subset [2^{\ell-3}, 2^{\ell}]$. 
Also, let $\beta\in \mathrm C_c^\infty((2^{-1}, 2))$ such that $\sum  \beta(2^{-k}\cdot) =1$ on $\mathbb R_{+}$. Then, we set 
\[   a_{\ell, j} (s,t,\xi) =     \tilde \chi_\ell (t)  \beta(2^{-j}|\xi|) \tilde a(s,t,\xi), \]  
so  $\sum_{\ell, j} a_{\ell, j}= \tilde a$. 
By  \eqref{fourierF},  $a_{\ell, j} =0$ if $|j-j_\ast|> C$ for some $C>0$. 

Denoting $(a)_\rho(s,t,\xi)=a(s,\rho t,\rho^{-1}\xi)$, via rescaling we observe  
$ \mathcal A_{\rho t}
[\gamma_{\sz}^{\delta}, \raisebox{-.2ex}{$a$} \,]g(x)= \mathcal A_{t}
[\gamma_{\sz}^{\delta}, \raisebox{-.2ex}{$(a)_{\rho}$} \,]g(\rho\, \cdot)( x/\rho).$
Thus, changes of variables yield
\[ \|  \mathcal A_{t}
[\gamma_{\sz}^{\delta}, \raisebox{-.2ex}{$a_{\ell, j}$} \,]\tilde f \|_{L^p(\mathbb R^{d+1})} = 2^{(\ell-2)/p} \| \mathcal A_{t}
[\gamma_{\sz}^{\delta}, \raisebox{-.2ex}{$(a_{\ell, j})_{2^{\ell-2}}$} \,] \tilde f_\ell  \|_{L^p(\mathbb R^{d+1})}, \] 
where $\tilde f_\ell= 2^{(\ell-2)d/p } \tilde f(2^{\ell-2}\cdot)$. 
Since $ \mathcal A_t  [\gamma_{\sz}^{\delta}, \raisebox{-.2ex}{$\tilde a$}]= \ssum{\ell, j} \mathcal A_t  [\gamma_{\sz}^{\delta}, \raisebox{-.2ex}{$a_{\ell, j}$}]$,  by \eqref{res:TA1} we get 
\[ \big\| \tilde\chi\, \mathcal T [\fa] f \big\|_{L^p(\mathbb R^{d+1})}\lesssim \delta  \ssum{\ell, j}\, \big \|   \mathcal A_t  [\gamma_{\sz}^{\delta}, \raisebox{-.2ex}{$(a_{\ell, j})_{2^{\ell-2}}$} \,] \tilde f_\ell  \big\|_{L^p(\mathbb R^{d+1})}. \]

To complete the proof, we only  have to relabel $(a_{\ell, j})_{2^{\ell-2}}$, $\ell=1,2,3$, $j_\ast-C\le j\le j_\ast +C$.  Indeed, since $\tilde a\in C^{d+N+2}$, $(a_{\ell, j})_{2^{\ell-2}}\in C^{d+N+2}$, which 
is supported in $I \times [2^{-1},4] \times \mathbb A_{j+\ell-2}$. Obviously,  \eqref{diff:a} holds for $\tilde a=(a_{\ell, j})_{2^{\ell-2}}$ because $\ell=1,2,3$. 
Changing variables $s\rightarrow \delta s+\sz$ and $\xi\to \delta^{-N}\widetilde{\mathcal L}_{\sz}^\delta\xi$ in \eqref{lowD}, by \eqref{tildegamma} we see that  \eqref{lowD} on $\supp_{s,\xi} \fa$ is equivalent to 
$
\sum_{j=1}^{N-1} | \langle (\gamma_{\sz}^\delta)^{(j)}(s), \xi \rangle|
\ge 2^k\delta^N/B'
$ for 
$(s,\xi) \in \supp_{s,\xi}   \mathfrak a_{\delta, \sz} ( \,\cdot\,,\delta^{-N}{\mathcal L}_{\sz}^\delta\,\cdot )$. 
Note   $\supp_{s,\xi}   \mathfrak a_{\delta, \sz} ( \,\cdot\,,\delta^{-N}{\mathcal L}_{\sz}^\delta\,\cdot )\supset \supp_{s,\xi}  \tilde a$. 
So, the same holds on $\supp_{s,\xi}  \tilde a$ and hence on  $\supp_{s,\xi}  (a_{\ell, j})_{2^{\ell-2}}$ if $B'$ replaced by $2B'$.   Therefore,  $ C^{-1}(a_{\ell, j})_{2^{\ell-2}}$ is of type $(j+\ell-2, N-1, \tilde B)$  relative to $\gamma_{\sz}^\delta$ for a large constant $C=C(B)$.  
\end{proof}

\subsection{Preliminary decomposition and reduction}
\label{sec:prop2.3} 
For the proof of Proposition \ref{newprop}, we make some reductions  by  decomposing the symbol $a$. 
We fix a sufficiently small positive  constant 
\[
\delta_*\le \min \{ \delta', (2^{7d}B^6)^{-N}\}, 
\]
which is to be specified in what follows. Here $\delta'$ is the number given in Lemma \ref{lem:res}.

We recall  that  $\gamma$ satisfies  \eqref{curveB},  $\ndeg N B$, $\vcon N B$,  and    $ a$ is of type $(k, N, B)$ relative to $\gamma$.  We set
\begin{align}
\label{a-N}
\eta_N(s,\xi)=\prod_{1\le j\le N-1} \beta_0 \Big( B2^{-k-1}\delta_*^{j-N} \langle\gamma^{(j)}(s),\, \xi \rangle \Big), 
\end{align}
where $\beta_0\in \mathrm C_c^\infty((-1,1))$ such that $\beta_0= 1$ on $[-1/2,1/2]$. 
It is easy to see 
$|\partial_s^{j} \partial^{l}_t \partial_{\xi}^{\alpha} (a\eta_N)| \le C |\xi|^{-|\alpha|}$
for $(j,l,\alpha)\in \mathcal I_N$, and the same holds for $a(1-\eta_N)$.

Note $\sum_{j=1}^{N-1}|\gamma^{(j)}(s)\cdot \xi|\ge (2B)^{-1} \delta_*^N|\xi|$ on  $\supp_{s,\xi}(a(1-\eta_N))$. So, we see   
$a(1-\eta_N)$ is a symbol of type $(k, N-1, B')$ for  $B'=CB^2 \delta_*^{-C}$ with a large $C$. 
Applying the assumption (Theorem \ref{lclb} with $L=N-1$ and $B=B'$), we  obtain 
\[ 
\| \mathcal A_t[\gamma, a(1-\eta_N)]f\|_{L^p(\mathbb R^{d+1})}\le
C 2^{(-\frac{2}p+\epsilon) k}\|f\|_{L^p(\mathbb R^d)}, \quad p\ge 4N-6.
\]
Thus,  it suffices  to consider $\mathcal A_t[\gamma,  a\eta_N]$.
Since $\ndeg{N}B$ holds on $\supp_{s,\xi} \fa$, 
\begin{align} | \gamma^{(N)}(s)\cdot \xi |\ge (2B)^{-1}|\xi|  \label{lowN}  \end{align}
whenever $(s,t,\xi)\in \supp a\eta_N$ for some $t$.

\medskip

\noindent\emph{Basic assumption} Before we continue  to prove the estimate for $\mathcal A_t[\gamma,  a\eta_N]$,  we make several assumptions which 
are clearly permissible by elementary decompositions. 

Decomposing $a$ we may assume that  $\supp_\xi a$  is contained in  a narrow conic neighborhood and $\supp_s a\subset I(\sz, \delta_\ast)$ for some $\sz$. 
  Let us set   
\[\Gamma_k =  \big\{\xi\in \mathbb A_k:   \dist\big( {|\xi|^{-1}}{\xi}, |\xi'|^{-1}{\xi'}\big)<  \delta_\ast \ \text{ for some } \ \xi'\in \supp_\xi (a\eta_N) \big\}. \]
We may also assume $\gamma^{(N-1)}(s') \cdot \xi'=0$ for some $(s', \xi')\in I(\sz, \delta_\ast)\times  \Gamma_k$. Otherwise,  $|\gamma^{(N-1)}(s) \cdot \xi|\gtrsim  |\xi|$ on 
$\supp_{s,\xi} a\eta_N$ and hence  $a\eta_N=0$ if we take $B$ large enough.  By  \eqref{lowN} and  the implicit function theorem,
there exists $\sigma $  such that
\begin{align}\label{zero}
 \gamma^{(N-1)}(\sigma(\xi)) \cdot \xi=0
   \end{align} 
 in a narrow conic  neighborhood of $\xi'$ where $\gamma\in C^{2d+2}$ since  $\gamma \in \mathrm C^{3d+1}(I)$. So,  decomposing $a$ further,  we may assume  
 $\sigma\in \mathrm C^{2d+2}(\Gamma_k)$ and $\sigma(\xi) \in  I(\sz, \delta_\ast) $ for $\xi \in \Gamma_k$. 
 Furthermore, since $\sigma$ is homogeneous of degree zero, 
\Be
\label{decaying-sigma}
 |\partial_\xi^\alpha \sigma(\xi)|  \le  C   |\xi|^{-|\alpha|}, \qquad \xi\in   \Gamma_k
\Ee 
for  a constant $C=C(B)$  if $|\alpha|\le 2d+2$.    Any symbol which appears in what follows is to be given by decomposing the symbol $a$ with appropriate cutoff functions. 
So, \emph{the $s,\xi$-supports of the symbols are assumed to be contained in   $I(\sz, \delta_\ast)\times\Gamma_k$.}
 
\medskip
We break $a$ to have further localization on the Fourier side. 
Let
\begin{align*}
\mathfrak a_1 (s,t,\tau,\xi)=
	a\eta_N \, \beta_0 \big( 2^{-2k}\delta_*^{-2N}|\tau+ \langle \gamma(s),\xi \rangle|^2\big)
	\end{align*}
and $\mathfrak a_0=a\eta_N-\mathfrak a_1$.
Then,  by  Fourier inversion
\[
  \mathcal A_t[\gamma,a\eta_N]f=\mathcal T[ \mathfrak a_1] f+ \mathcal T [ \mathfrak a_0]f.
\]
It is easy to show  $ \| \mathcal T[ \mathfrak a_0]f \|_{p} 
\lesssim_B\! 2^{-2k} \|f\|_{p}$ for $1 \le p \le \infty.$
Indeed, consider 
$  \tilde {\mathfrak a}_0 =  -(\tau+\gamma(s)\cdot \xi)^{-2} {\partial_{t}^2\mathfrak a_0}. 
$
By \eqref{Ma} and integration by parts in $t'$,   $m[\mathfrak a_0]=m[\tilde {\mathfrak a}_0]$ and hence  $\mathcal T [ \mathfrak a_0 ]= \mathcal T [\tilde {\mathfrak a}_0]$.  
Thanks to \eqref{conv},  it is sufficient to show 
\begin{align*} 
| K [ \tilde{\mathfrak a}_0 ](s,t,x) \big| 
	&\le C\, 
2^{k(d-1)}  \!
\int  \! \big(1+2^k  |t-t'|  +2^k 
 |x-t'\gamma(s)|\big)^{-d-3}\,dt'
\end{align*}
for a constant $C=C(B, \delta_\ast)$. 
Note $|\tau+\langle \gamma(s),\xi \rangle|\gtrsim  2^k$ on $\supp \tilde {\mathfrak a}_0$, and recall
 \eqref{eq:kernel}.  Rescaling and integration by parts in $\tau,\xi$, as in the proof of Lemma \ref{kernel}, show the estimate.

The difficult part is to estimate    $\mathcal T[\mathfrak a_1]$. 
Since $\delta_*$ is a fixed constant, it is obvious that $C^{-1}\mathfrak a_1 \in \mathfrak A_k(\sz, \delta_*)$ for some $C=C(B,\delta_\ast)$. 
 So, the desired  estimate for $\mathcal T[\mathfrak a_1]$ follows once we have  the next proposition.

\begin{prop}\label{microlocal}
Let  $ \mathfrak a\in \mathfrak A_k(\sz, \delta_*)$ with $\supp_{\xi} \fa \subset  \Gamma_k$. 
Suppose  Theorem \ref{lclb} holds for $L=N-1$. 
Then, if   $p\ge 4N-2$, for $\epsilon>0$ we have
\begin{align}\label{overline-a0}
\big\| \mathcal T[\mathfrak a]f \big\|_{L^p(\mathbb R^{d+1})} 
\le C_\epsilon 2^{-\frac{2} pk+ \epsilon k}\|f\|_{L^p(\mathbb R^d)}.
\end{align}
\end{prop}

Therefore, the proof of Proposition \ref{newprop} is completed if we prove Proposition \ref{microlocal}. 
For the purpose, we use  Proposition \ref{iterative} below,  which allows us  to decompose $ \mathcal T[\mathfrak a]$  into operators given by symbols with smaller $s$-supports while the consequent minor parts have acceptable bounds. 
This type of argument  was used in \cite{PS} 
when $L=2$.

Let  $\delta_0$ and $\delta_1$ be positive numbers such that
\Be\label{orderB}
2^{7d}B^6\delta_0^{(N+1)/N}  \le \delta_1\le \delta_0 \le \delta_*,   \qquad 2^{-k/N}\le \delta_1.  
\Ee
Then, it is clear that 
\Be
\label{ord2} 
B^{6N}\delta_0^{j+1} \le 2^{-7dN} \delta_1^{j}, \qquad j=1,\dots, N.
\Ee

For $n\ge 0$, we denote  $
\mathfrak J_n^\mu=\{ \nu \in \mathbb Z: |2^n\delta_1 \nu - \delta_0 \mu| \le \delta_0\}.
$

\begin{prop}\label{iterative}
For $\mu$ such that  $\delta_0\mu \in I(\sz, \delta_\ast)\cap \delta_0\mathbb Z$, 
let $\mathfrak a^\mu \in \mathfrak A_k(\delta_0\mu,\delta_0)$ with $\supp_{s,\xi} \fa^\mu\subset I(\sz, \delta_\ast)\times\Gamma_k$. Suppose Theorem \ref{lclb} holds for $L=N-1$. 
 Then, if  $p\ge 4N-2$, for $\epsilon>0$
there exist a constant $C_\epsilon=C_\epsilon (B) \ge2$ and symbols $\mathfrak a_{\nu} \in \mathfrak A_k(\delta_1\nu,\delta_1)$ with $\supp_{s,\xi} \fa_\nu\subset I(\sz, \delta_\ast)\times\Gamma_k$,
$\nu \in \cup_\mu\mathfrak J_0^\mu$, such that  
\begin{align*}
\!\!\big( \sum_\mu \| \mathcal T[ \mathfrak a^\mu]f\|_p^p\big)^{\frac1p}\!
\le C_\epsilon
	& \big({\delta_1}/{\delta_0} \big)^{\frac {2N}{p}-1-\epsilon}
	 \big( \sum_\nu \| \mathcal T[\mathfrak a_{\nu}]f\|_p^p \big)^{\frac1p}
\!+C_\epsilon \delta_0^{-\frac {2N}p +1+\epsilon } 2^{-\frac{2}pk+2\epsilon k}\|f\|_p. 
\end{align*}
\end{prop}

Assuming Proposition \ref{iterative}, we prove Proposition \ref{microlocal}.

\subsection{Proof of Proposition \ref{microlocal}}
\label{sec:prop2.5}
Let $\mathfrak a \in \mathfrak A_k(\sz,\delta_*)$.
We may assume ${\sz}=\delta_\ast\mu$ for some $\mu \in \mathbb Z$.
To apply  Proposition \ref{iterative} iteratively, we need to choose 
an appropriate decreasing sequence of positive numbers since the decomposition is subject to the condition  \eqref{orderB}. 

Let $\delta_0=\delta_*$, so $(2^{7d}B^6)^{N}\delta_0< 1$. 
 Let $J$ be the largest integer such that  
\[
(2^{7d}B^6)^{N(\frac{N+1}N)^{J-1}-N}  \delta_0^{(\frac{N+1}N)^{J-1}} >2^{-\frac kN}.
\]
So, $ J \le C_1 \log k$ for a constant $C_1\ge1$.  We  set  
\begin{align}
\label{defn-j}
 \delta_{\!J}=2^{-\frac kN}, \qquad  \delta_j=(2^{7d}B^6)^{N(\frac{N+1}N)^{j}-N} &\delta_0^{(\frac{N+1}N)^{j}} 
\end{align}
for $j=1,\dots,J-1$. Thus, it follows that  
\begin{align} \label{del-con} 
&2^{7d}B^6\delta_j^{(N+1)/N}\le \delta_{j+1} <\delta_j,   \qquad \, j=0, \dots, J-1.
\end{align}

For a given $\epsilon>0$, let $\tilde \epsilon=\epsilon/4$.
Since $\mathfrak a\in \mathfrak A_k( \delta_0\mu,\delta_0)$  and \eqref{orderB} holds for  $\delta_0$ and $\delta_1$,  applying Proposition \ref{iterative} to $\mathcal T[\mathfrak a]$, we have 
\begin{align*}
 \| \mathcal T[\mathfrak a]f\|_p
\le C_{\tilde \epsilon} 
	& \big({\delta_1}/{\delta_0}\big)^{\frac {2N}{p}-1-\tilde \epsilon}
	 \big(\sum_{\nu_1} \| \mathcal T[ \mathfrak a_{\nu_1}]f\|_p^p\,\big)^{\frac1p}
+C_{\tilde \epsilon}
	\delta_0^{-\frac {2N}p+1 +\tilde \epsilon } 2^{-\frac{2}pk+2\tilde \epsilon k}\|f\|_p,
\end{align*}
where $ \mathfrak a_{\nu_1} \in \mathfrak A_k(\delta_1\nu_1,\delta_1)$, $\nu_1 \in \mathfrak J_0^\mu$.  
Thanks to \eqref{del-con}  we  may apply  again Proposition \ref{iterative} to $\mathcal T[\mathfrak a_{\nu_1}]$ while  $\delta_0$, $\delta_1$ replaced by 
$\delta_1,$ $\delta_2$, respectively.  Repeating this procedure up to $J$-th step yields symbols $\mathfrak a_\nu \in \mathfrak A_k(\delta_{\!J}\nu, \delta_{\!J})$, $\delta_{\!J} \nu \in \delta_{\!J}\mathbb Z\cap I(\delta_0\mu, \delta_0)$,  such that 
\[ 
\| \mathcal  T[\fa]f\|_p
	\le C_{\tilde \epsilon}^J \delta_{\!J}^{\frac {2N}{p} -1-\tilde \epsilon}
	 \big(\sum_\nu  \|\mathcal T[\mathfrak a_\nu]f\|_p^p \,\big)^{\frac1p}
	+\sum_{0\le j\le J-1} C_{\tilde \epsilon}^{j+1} \delta_{0}^{-\frac{2N}p+1+\tilde \epsilon} 2^{-\frac{2} pk+ 2\tilde \epsilon k} \|f\|_p
		\]
for $p \ge 4N-2$.  Now, assuming  
\begin{equation}\label{est03}
 \big( \sum_\nu\| \mathcal T[ \mathfrak a_{\nu} ]f\|_p^p \,\big)^{1/p}
	\lesssim_B  	2^{-k/N} 	\|f\|_p, \qquad  2\le p\le \infty 
\end{equation}
for the moment, we can finish the proof of Proposition \ref{microlocal}.
Since $C_{\tilde \epsilon}\ge 2$, combining the above inequalities, 
we get
\[ 
\| \mathcal  T[\mathfrak a]f\|_p
\lesssim_B \!C_{\tilde \epsilon}^{J+1} 
\big( 2^{-\frac 2p k+\frac{\tilde \epsilon }{N}k}
+2^{-\frac{2} pk +2\tilde \epsilon k}\big) \|f\|_p.
\]
Note $J \le C_1 \log k$, so 
$C_{\tilde \epsilon}^{J+1}  \le C' 2^{\epsilon k/2}$  for some $C'$ if $k$ is sufficiently large.
 Thus, the right hand side is bounded  by 
$C 2^{- 2k/p+\epsilon k} \|f\|_p$.

\medskip

It remains to show  \eqref{est03} for $2\le p\le \infty$. By  interpolation, it is enough to obtain   
 \eqref{est03} for $p=\infty$ and $p=2$. 
The case $p=\infty$ follows by \eqref{ker-est}   since  $\mathfrak a_\nu \in \mathfrak A_k(\delta_{\!J}\nu,\delta_{\!J})$. So, we need only to prove  \eqref{est03}
 for $p=2$. To do this,  we first observe the following,  which shows 
 $\supp_\xi \fa_{\nu}$ are finitely overlapping.

\begin{lem} 
\label{disjoint}
For $b\ge 1$, $s\in I(\sz, \delta_\ast)$,  and  $0<\delta \le \delta_*$, let us  set
\begin{align}
\label{wl}
{\Lambda}'_k(s,\delta, b) = \bigcap_{1\le j\le N-1} \big\{ \xi  \in   \Gamma_k
 : | \langle \gamma^{(j)}(s), \xi \rangle| \le b 2^{k} \delta^{N-j} \big \}.
\end{align}
If  ${\Lambda}'_k(s_1, \delta, b) \cap {\Lambda}'_k(s_2, \delta, b)\neq \emptyset$ for some $s_1, s_2\in  I(\sz, \delta_\ast)$, then there is a constant $C=C(B)$ such that   $|s_1-s_2|\le  Cb\delta $. 
\end{lem}

\begin{proof}  
Let $\xi \in  \Lambda_k'(s_1,\delta,b) \cap  \Lambda_k'(s_2,\delta,b)$. 
Since $|\gamma^{(N-1)}(s_j) \cdot \xi |\le b 2^k\delta$, $j=1,2$,   by \eqref{zero} and \eqref{lowN} we see $|s_j-\sigma(\xi)|\le   2^2bB \delta$, $j=1,2$, using the mean value theorem.  
This implies $|s_1-s_2|\le    2^3bB \delta $.  
\end{proof}

We recall \eqref{Ma}. Since \eqref{lowN} holds on $\supp_{s,\xi} \fa_\nu$,  by van der Corput's lemma (e.g., see \cite[ Corollary, p. 334]{Stein93})
we have 
\begin{align*}
|m[\mathfrak a_{\nu}](\tau,\xi)| \lesssim 2^{-k/N} 
\big( \| \mathfrak a_{\nu} (\cdot,t,\tau,\xi) \|_\infty
+\| \partial_s \mathfrak a_{\nu} (\cdot, t,\tau,\xi) \|_1\big)
\lesssim_B\! 2^{-k/N}.
\end{align*}
The second inequality is clear since $\mathfrak a_{\nu} \in \mathfrak A_k(\delta_{\!J}\nu, \delta_{\!J})$.  From \eqref{Ta} note $\mathcal F(\mathcal T[\mathfrak a_{\nu}] f)=  m[\mathfrak a_{\nu}]\widehat f$.  Since $\supp_\xi  \mathfrak a_{\nu}\subset S_\nu:= {\Lambda}'_k(\delta_{\!J}\nu, \delta_{\!J}, 2^5 B)$, 
 $\supp_\xi  \mathcal F(\mathcal T[\mathfrak a_{\nu}] f)\subset S_\nu$  by \eqref{Ma}. By Lemma \ref{disjoint} it follows that 
the sets $S_\nu$ overlap at most $C=C(B)$ times.   Therefore, 
 Plancherel's theorem and the estimate above yield   
\[  \| \sum_\nu \mathcal T[\mathfrak a_{\nu}] f\|_2^2
\lesssim_B\!  2^{-2k/N}   \sum_\nu \int_{S_\nu}\int_{\{\tau: | \tau  + \gamma(\delta_{\!J}\nu)\cdot \xi|\le 2^5 B\}} \,d\tau\, |\widehat f(\xi)|^2 \,d\xi
\]
since   $\supp \mathfrak a_{\nu}\subset \Lambda_k(\delta_{\!J}\nu, \delta_{\!J}, B)$.   This gives  \eqref{est03} for $p=2$.

\subsection{Decoupling inequalities}
\label{sec:dec} 
We denote   $\mathbf r_\circ^N(s)=(s,s^2/2!,\,\dots\, s^N/N!)$, and consider  a collection of curves from $I$  to $\mathbb R^N$ which are small perturbations of  $\mathbf r_\circ^N$:    
\[\mathfrak C(\epsilon_\circ;N) :=\{ \mathbf r \in \mathrm C^{2N+1}(I):
\|\mathbf r-\mathbf r_\circ^N \|_{\mathrm C^{2N+1}(I)} < \epsilon_\circ\}. \] 

For $\mathbf r\in \mathfrak C(\epsilon_\circ;N) $ and $s\in I$, we define 
\[ \mathcal N_\mathbf r(s,\delta)=  \Big\{ \mathbf r(s)+ \sum_{1\le j\le N}  u_j { \mathbf r^{(j)}(s)}: \ |u_j|\le \delta^j, \quad~j=1,\dots,N \Big\}. \]  
 Let $s_1, \dots, s_l\in I$ be $\delta$-separated points, i.e.,  $|s_n-s_j|\ge \delta$ if $n\neq j$, such that   $\bigcup_{j=1}^l(s_j-\delta, s_j+\delta)\supset I$. Then, we set 
 \[ \theta_j=\mathcal N_\mathbf r(s_j,\delta), \qquad  1\le j\le l.  \]
 
The following is due to  Bourgain, Demeter, and Guth \cite{BDG} (also see   \cite{GLYZ}).

\begin{thm}
Let $0<\delta\ll 1$.   Suppose  $\mathbf r\in \mathfrak C(\epsilon_\circ;N) $ for a small enough $\epsilon_\circ>0$. 
Then, if   $2 \le p\le N(N+1)$,   for $\epsilon>0$ we have
\begin{align}
\label{eq:decoupling}
\big\| \sum_{1\le j\le l}  f_j \big\|_{L^p(\mathbb R^{N})}\le C_\epsilon  \delta^{-\epsilon}
\big( \sum_{1\le j\le l}\|f_j \|_{L^p(\mathbb R^{N})}^2
\big)^{1/2}
\end{align}
whenever  $\supp \widehat f_j \subset \theta_j$ for $1\le j\le l$. 
\end{thm}

The constant $C_\epsilon$ can be taken to be independent of particular choices of the $\delta$-separated points $s_1, \dots, s_l$. 
One can obtain a conical extension of the inequality \eqref{eq:decoupling} by modifying  
the argument in \cite{BD} which deduces the decoupling inequality for the cone from that for  the paraboloid (see \cite[Proposition 7.7]{BGHS2}). 
Let us consider  conical sets  
 \[\otheta_j=\{ (\eta, \rho)\in \mathbb R^N\times [1,2]:  \eta/\rho\in \theta_j\}, \quad 1\le j\le l.\]

\begin{cor}\label{cdecoupling}
Let $0<\delta\le1$ and let $\mathbf r\in \mathfrak C(\epsilon_\circ;N) $ with a small enough $\epsilon_\circ>0$.  
Then, if  $2 \le p\le N(N+1)$,  for $\epsilon>0$  we have 
\begin{align}
\label{c-decoupling}
\big\|\sum_{1\le j\le l}  F_j \big\|_{L^p(\mathbb R^{N+1})}\le C_\epsilon \delta^{-\epsilon}
\big(\sum_{1\le j\le l} \| F_j\|_{L^p(\mathbb R^{N+1})}^2
\,\big)^{1/2}
\end{align} 
whenever   
$\supp \widehat F_j \subset \otheta_j$ for $1\le j\le l$. 
\end{cor}

The decoupling inequality \eqref{c-decoupling} does not fit  the symbols which appear later when we decompose $\mathfrak a$ (see  Section \ref{sec:decomp} and Section \ref{pf:thm1}).
As to be seen later,  those symbols are  related to the slabs of the following form.

\begin{defn} 
Let $N\ge2$ and $\tilde{\mathbf r}\in \mathfrak C(\epsilon_\circ;N+1)$. 
For $s\in I$, we denote by  $\mathbf s(s,\delta, \rho ;\tilde{\mathbf r})$ 
 the set of $(\tau,\eta) \in  \mathbb R\times\mathbb R^{N}$ which satisfies   
\begin{align*}
 & \rho^{-1}\le|\langle \tilde{\mathbf r}^{(N+1)}(s), (\tau,\eta) \rangle|\le 2\rho, 
 \\[3pt]
 &\quad  |\langle \tilde{\mathbf r}^{(j)}(s), (\tau,\eta) \rangle| \le \delta^{N+1-j}, \qquad  \  \  j=N,\dots,1, 
\end{align*}
\end{defn}

The same form of  decoupling  inequality continues to be  valid for the slabs $\mathbf s(s_1,\delta,1; \tilde{\mathbf r}), \dots,$ $\mathbf s(s_l,\delta,1; \tilde{\mathbf r})$.
Beltran et al. \cite[Theroem 4.4]{BGHS2} showed, using the Frenet–Serret formulas,  that 
those slabs can be  generated by conical extensions of the slabs given by a nondegenerate curve in $\mathbb R^N$.  The following is a consequence of 
Corollary \ref{cdecoupling} and a simple manipulation using decomposition and rescaling. 

\begin{cor}\label{rcdecoupling}
Let $0<\delta\le1$, $\rho\ge 1$,  and  $\tilde{\mathbf r} \in \mathfrak C(\epsilon_\circ;N+1)$ for a small enough $\epsilon_\circ>0$. 
Denote $\mathbf s_j=\mathbf s(s_j,\delta, \rho ;\tilde{\mathbf r})$ for $1\le j\le l$. 
 Then, if   $2 \le p\le N(N+1)$,  for $\epsilon>0$ there is a constant $C_\epsilon=C_\epsilon(\rho)$ such that 
 \begin{align}
\label{c-decoupling2}
\big\| \sum_{1\le j\le l}  F_j \big\|_{L^p(\mathbb R^{N+1})}\le C_\epsilon  \delta^{-\epsilon}
\big( \sum_{1\le j\le l} \|F_j \|_{L^p(\mathbb R^{N+1})}^2
\big)^{1/2}
\end{align} 
whenever  
$\supp \widehat {F_j} \subset \mathbf s_j$ for $1\le j\le l$.
\end{cor}

For our purpose of proving Proposition \ref{iterative},  we use a modified form. 
If $p_\ast\in  [2 , N(N+1)] $, then we have
\[ 
\big\| \sum_{1\le j\le l}  F_j \big\|_{L^p(\mathbb R^{N+1})}\le C_\epsilon 
  \delta^{-1+\frac{2+p_\ast}{2p}-\epsilon}
\big( \sum_{1\le j\le l} \|F_j \|_{L^p(\mathbb R^{N+1})}^p
\big)^{1/p}
\]
for $p\ge p_\ast$. The case $p=p_\ast$ follows by \eqref{c-decoupling2} and H\"older's inequality. Interpolation  with the trivial $\ell^\infty L^\infty$--$L^\infty$  estimate gives  the estimate  for $p> p_\ast$.  One may  choose different $p_\ast$   for the particular purposes.  In fact,  for the local smoothing estimate we take $p_\ast= 4N-2$ to have 
\begin{align}
\label{dgain0}
\big\| \sum_{1\le j\le l}  F_j \big\|_{L^p(\mathbb R^{N+1})}\le  C_\epsilon    \delta^{-1+\frac{2N}{p}-\epsilon}
\big( \sum_{1\le j\le l} \|F_j \|_{L^p(\mathbb R^{N+1})}^p
\big)^{1/p}
\end{align} 
for  $p\ge 4N-2$
(see Section \ref{sec:pf-4.1}).  For the $L^p$ Sobolev regularity estimate, we observe 
\begin{align}
\label{dgain}
\big\| \sum_{1\le j\le l}  F_j \big\|_{L^p(\mathbb R^{N+1})}\le  C_{\epsilon_0}   \delta^{-1+\frac{N+1}{p}+\epsilon_0}
\big( \sum_{1\le j\le l} \|F_j \|_{L^p(\mathbb R^{N+1})}^p
\big)^{1/p}
\end{align} 
holds for some $\epsilon_0= \epsilon_0(p)>0$ if $2N<p<\infty$. Indeed, we need only  to take $p_\ast> 2N$ close enough to $2N$.   
The presence of $\epsilon_0$ in  \eqref{dgain}  is crucial for proving the optimal  Sobolev regularity estimate (see Proposition \ref{iterative2}). 

The inequalities \eqref{dgain0} and \eqref{dgain}  obviously  extend to  cylindrical forms via the Minkowski inequality. For example, 
set $  \tilde{\mathbf s}_j=\big\{ (\xi,\eta )\in \mathbb R^{N+1}\times \mathbb R^M:  \xi\in \mathbf s_j  \big\}$   
for  $1\le j\le l$. 
Then, using \eqref{dgain0},  we have 
\Be
\label{cylinder}
 \big\| \sum_{1\le j\le l}  G_j \big\|_{L^p(\mathbb R^{N+M+1})}\le C_\epsilon    \delta^{-1+\frac{2N}{p}-\epsilon}
\big( \sum_{1\le j\le l}  \| G_j \|_{L^p(\mathbb R^{N+M+1})}^2
\big)^{1/2} \Ee
whenever  
$\widehat {G}_j$ is supported 
in $\tilde{\mathbf s}_j$. Clearly, we also have a similar extension of \eqref{dgain}. 

\section{Decomposition of  the symbols}\label{sec:mainthm}
In this section, we prove Proposition \ref{iterative}  by applying the decoupling inequality. Meanwhile, the induction assumption (Theorem \ref{lclb} with $L=N-1$)  plays an important role.  We decompose a given symbol  $\mathfrak a^\mu \in \mathfrak A_k(\delta_0\mu,\delta_0)$   
into  the symbols  with their $s$-supports contained in  intervals of length about $\delta_1$ while  the consequent  minor contribution  is controlled within an acceptable bound.   
To achieve it up to $\delta_1$ satisfying \eqref{orderB}, 
we   approximate  $\langle G(s), (\tau,\xi) \rangle $ in a local coordinate system near the set $\{(s,\xi):\langle \gamma^{(N-1)}(s),\xi\rangle=0\}$. 

\subsection{Decomposition of the symbol $\mathfrak a^\mu$} 
\label{sec:decomp}
We begin by introducing some notations. 

Fixing $\mu \in \mathbb Z$ such that $\delta_0\mu\in I(\sz, \delta_\ast)$, we  consider  linear maps 
\[
y_\mu^j(\tau,\xi)= \langle G^{(j)}(\delta_0\mu ),(\tau,\xi)\rangle,  \qquad j=0, 1,\dots,N.
\]
In particular,  $y_\mu^j(\tau,\xi)=\langle\gamma^{(j)}(\delta_0\mu),\xi\rangle$ if $1\le j\le N.$  By \eqref{lowN} it follows that 
\begin{align}
\label{YN}
|y_\mu^N(\tau,\xi)| \ge (2B)^{-1}|\xi| .
\end{align}

We  denote  
\[
\omega_\mu(\xi)=\frac{y_\mu^{N-1}(\tau,\xi)}{y_\mu^{N}(\tau,\xi)}, 
\]
which is  close to  $\delta_0\mu-\sigma(\xi)$ (see \eqref{s}).  Then, we define $\mathfrak g_\mu^N, \mathfrak g_\mu^{N-1},\dots, \mathfrak g_\mu^{0}$  recursively, by setting 
$
\mathfrak g_\mu^N=y_\mu^N,$ and 
\Be
\label{def-g}
\mathfrak g_\mu^{j}(\tau,\xi)=y^{j}_\mu(\tau,\xi)-\sum_{\ell=j+1}^{N}\frac{\mathfrak g^{\ell}_\mu(\tau,\xi)}{(\ell-j)!}
(\omega_\mu(\xi))^{\ell-j},  \quad j=N-1, \dots, 0. 
\Ee
Note that $\mathfrak g_\mu^{N-1}=0$ and \eqref{def-g} can be rewritten as  follows:
\begin{align}\label{formulaw}
y^{m}_\mu(\tau,\xi)
=\sum_{\ell=m}^{N}\frac{\mathfrak g^{\ell}_\mu(\tau,\xi)}{(\ell-m)!}
(\omega_\mu(\xi))^{\ell-m},
\quad m=0,\dots,N.
\end{align}
The identity continues to  hold for $m=N$ since $\mathfrak g_\mu^N=y_\mu^N$.   
Apparently, $\mathfrak g_\mu^1,\dots, \mathfrak g_\mu^N$ are independent of $\tau$ 
since so are $y_\mu^1, \dots, y_\mu^N$.

 \medskip
 
For $j=1,\dots,N$, set
\Be 
\label{ej} 
\mathcal E_j(\xi):=(y_\mu^N(\tau,\xi))^{-1}\int_{\sigma(\xi)}^{\delta_0\mu}\frac{\langle\gamma^{(N+1)}(r),\xi\rangle}{j!}(\sigma(\xi) -r)^{j}dr.
\Ee
By \eqref{ej} with $j=1$ and  integration by parts, we have
\Be
\label{s} \mathcal E_1(\xi) = \sigma(\xi) -\delta_0\mu+\omega_\mu(\xi). 
\Ee

\begin{lem} 
For $0\le j\le N-1$, we have   
\begin{align}
\label{uj}
	\langle G^{(j)}(\sigma(\xi)),(\tau,&\xi) \rangle =  \sum_{\ell=j}^{N}\frac{\,\, \mathfrak g^{\ell}_\mu\,\,(\mathcal E_1)^{\ell-j}}{(\ell-j)!} 
	 -y^N_\mu\mathcal E_{N-j}.
	\end{align}
\end{lem}

\begin{proof} 
When $j=N-1$, \eqref{uj}  is clear.  To show \eqref{uj} for $j=0,1,\dots,N-2$,  by Taylor's theorem with integral remainder  we have 
\[\langle G^{(j)}(\sigma(\xi)),(\tau,\xi)\rangle =\sum_{m=j}^{N}y^{m}_\mu (\tau,\xi) \frac{(\sigma(\xi)-\delta_0\mu)^{m-j}}{(m-j)!}-y^N_\mu(\tau,\xi)\mathcal E_{N-j}(\xi).\]
Using \eqref{formulaw} 
and then changing the order of the sums,  we see
\[
\langle G^{(j)}(\sigma(\xi)),(\tau,\xi)\rangle= \sum_{\ell=j}^{N} \mathfrak g^{\ell}_\mu
 \Big(\sum_{m=j}^{\ell}
 \frac{(\sigma(\xi)-\delta_0\mu)^{m-j}}{(\ell-m)!(m-j)!}
(\omega_\mu)^{\ell-m} \Big)
-y^N_\mu\,\mathcal E_{N-j}.
\]
The sum over $m$  equals 
$(\sigma(\xi)-\delta_0\mu+ \omega_\mu)^{\ell-j}
/(\ell-j)!$. So, 
   \eqref{uj} follows by  \eqref{s}. 
\end{proof}

We now decompose the symbol $\mathfrak a^\mu \in \mathfrak A_k(\delta_0\mu,\delta_0)$ by making use of  $\mathfrak g_\mu^j$, $j=0,\dots,N-2$. We  define  
\Be
\label{gn}
\mathfrak G_N^\mu(s,\tau,\xi)=\sum_{j=0}^{N-2}\big(2^{-k} \mathfrak g^{j}_\mu(\tau,\xi)\big)^{\frac{2N!}{N-j}}+(s-\sigma(\xi))^{2N!}.
\Ee
Let $ \beta_N=\beta_0-\beta_0(2^{2N!}\cdot)$, so $\sum_{\ell \in \mathbb Z}\beta_N(2^{2N!\ell}\cdot )=1$. We also take $\zeta \in \mathrm C_c^\infty((-1,1))$ such that $\sum_{\nu \in \mathbb Z} \zeta(\cdot-\nu)=1$. 
For $n\ge 0$ and $\nu \in \mathfrak J_n^\mu$,   we set 
\[
\mathfrak a_{\nu}^{\mu,n}
= \mathfrak a^\mu 
\times
\begin{cases}
	\beta_0 \big( 
	\delta_1^{-2N!} \, \mathfrak G_N^\mu  
	\big) \, \zeta(\delta_1^{-1}s-\nu),   \quad &  n=0,
	\\[4pt]
	\beta_N \big( (2^{n}\delta_1)^{-2N!}
	\, \mathfrak G_N^\mu 
	\big) \,	\zeta(2^{-n}\delta_1^{-1}s-\nu),  \quad  & n\ge 1. 
	\end{cases}
\] 
  Then, it follows that 
\begin{align}\label{split}
 \mathfrak a^\mu = \sum_{n\ge 0} \sum_{\nu \in \mathfrak J_n^\mu} \mathfrak a^{\mu,n}_{\nu}.
\end{align} 
 
 \begin{lem}\label{itsym}  
There is a constant $C=C(B)$ such that $C^{-1}\mathfrak a_{\nu}^{\mu,n}\in \mathfrak A_k(2^n\delta_1\nu, 2^n\delta_1)$ for $n\ge 0$, $\mu,$ and $\nu$. 
\end{lem}
The proof of Lemma \ref{itsym} is elementary though it is somewhat  involved. We postpone the proof until Section \ref{pf-lem}.

We  collect some elementary  facts regarding $\mathfrak a_{\nu}^{\mu,n}$. First, we may assume 
\Be
\label{nnn}
2^n\delta_1\lesssim_B\!\delta_0
\Ee
 since, otherwise, $\mathfrak a_{\nu}^{\mu,n}=0$.     Note $|\langle \gamma^{(N-1)}(\delta_0\mu ),\xi\rangle|\le B2^{k+5}\delta_0$ if  $\xi\in \supp_{\xi} \mathfrak a_\mu$. 
  Then, \eqref{lowN}, \eqref{zero}, and the mean value theorem show  
  \Be 
 |\sigma(\xi)-\delta_0\mu|\le  B^22^{7}\delta_0 
  \label{sigma}
\Ee
for $\xi\in \supp_{\xi} \mathfrak a_\mu.$  
  If $(\tau, \xi)\in \supp_{\tau, \xi} \mathfrak a_\mu\subset \Lambda_k(\delta_0\mu, \delta_0, B)$, 
$|y_\mu^j(\tau,\xi)|\le B2^{k+5}\delta_0^{N-j} $ for $0\le j\le N-1$. Since $|\omega_\mu|\lesssim B^2\delta_0$, \eqref{def-g} gives $|\mathfrak g_\mu^j(\tau,\xi)|\lesssim_B\! 2^{k+5}\delta_0^{N-j} $ for $0\le j\le N-2$.  
Therefore,  
$ \mathfrak G_N^\mu\lesssim_{B}\!\delta_0^{2N!}$ on the support of $\mathfrak a_\mu$. This gives \eqref{nnn}. 

Since $\mathfrak G_N^\mu \le (2^n\delta_1)^{2N!}$ on $\supp \mathfrak a_{\nu}^{\mu,n}$, the following hold on the support of $\mathfrak a_{\nu}^{\mu,n}$:
\begin{align}
|s-\sigma(&\xi)|\le 2^n\delta_1, 
\label{gbound2}
\\[2pt]
\label{gbound}
2^{-k}|\mathfrak g^{j}_\mu(\tau,\xi)&|\le (2^n\delta_1)^{N-j},   \qquad  0\le j \le N-1.
\end{align} 
The inequality \eqref{gbound}  holds true for $j=N-1$ since  $\mathfrak g^{N-1}_\mu=0$.  We also have 
 \begin{align}
  \label{ejb} &|\mathcal E_j(\xi)|\le B^2 (B^22^7\delta_0)^{j+1}, 
  \\
  \label{si-2n} 
&|\sigma(\xi)-2^n\delta_1\nu| \le 2^{n+1}\delta_1.
\end{align}
on $\supp_\xi \mathfrak a_{\nu}^{\mu,n}$.    By using \eqref{ej}, \eqref{sigma}, and \eqref{YN}, it is easy to show \eqref{ejb}.  
Since  $|s-2^n\delta_1\nu| \le 2^n\delta_1$ on $\supp_s \mathfrak a_{\nu}^{\mu,n}$, 
 \eqref{si-2n}  follows by \eqref{gbound2}.

\subsection{Proof of Proposition \ref{iterative}}   
\label{sec:pf-4.1}
By \eqref{split} and the Minkowski inequality we have
\begin{equation}
\begin{aligned}\label{apnorm}
 \big( \sum_\mu \big\| \mathcal T [\mathfrak a^\mu]f \big\|_p^p \,\big)^{1/p}
	\le
\sum_{n\ge0}  \big(\sum_{\mu} \big\|   \sum_{\nu \in \mathfrak J_n^\mu} \mathcal T [ \mathfrak a^{\mu,n}_{\nu}]f \big\|_p^p\,\big)^{1/p}.
\end{aligned}
\end{equation}
We apply the inequality \eqref{dgain0} 
to $\sum_{\nu \in \mathfrak J_n^\mu} \mathcal T  \mathlarger[\mathfrak a^{\mu,n}_{\nu} \mathlarger ]f$ after a suitable linear change of variables. 
The  symbols $\mathfrak a_{\nu}^{\mu,0}$  are to constitute the set  $\{\mathfrak a_\nu\}$ while the operators associated to $\mathfrak a_{\nu}^{\mu,n}$, $n\ge 1$  are to be handled similarly as  in  Section 2. 

\subsubsection*{Applying the decoupling inequality} To prove Proposition \ref{iterative},
we  first show 
\begin{align}
\label{decoupn}
\big\| \sum_{\nu\in \mathfrak J_{n}^\mu} \mathcal T [\mathfrak a_{\nu}^{\mu,n}]f \big\|_p
	\le C_\epsilon 
	\big({2^n\delta_1}/{\delta_0} \big)^{\frac {2N} {p}-1-\epsilon}  \big(\sum_{\nu\in \mathfrak J_{n}^\mu} 
	\big\| \mathcal T [\mathfrak a_{\nu}^{\mu,n}]f \big\|_p^p \,\big)^{1/p}
\end{align}
for $p \ge 4N-2$.  To apply  the inequality   \eqref{dgain0}, 
we consider $\supp_{\tau,\xi}  \mathfrak a_{\nu}^{\mu,n}$, which contains 
the Fourier support  of $ \mathcal T [\mathfrak a_{\nu}^{\mu,n}]f$ as is clear from \eqref{Ma} and \eqref{Ta}.

We set
\[
\mathbf y_\mu(\tau,\xi)=\big(y_\mu^0(\tau,\xi),\dots,y^{N}_\mu(\tau,\xi)\big).
\]

\begin{lem} 
\label{supp-a}
Let $\mathbf r
 =\mathbf r_{\circ}^{N+1}$ and $\mathcal D_\delta$ denote the matrix  $ (\delta^{-N} e_1, \delta^{1-N} e_2, \dots,  \delta^0 e_{N+1})$ where $e_j$ denotes the $j$-th standard  unit vector in $\mathbb R^{N+1}$.    On $\supp_{\tau,\xi} \mathfrak a^{\mu,n}_{\nu}$,  we have 
     \begin{align}
     \label{yy0}
\Big|\Big\langle   \mathcal D_{\delta_0} \mathbf y_\mu(\tau,\xi), \mathbf r^{(j)}\Big(\frac{2^n\delta_1}{\delta_0}\nu-\mu\Big)\Big\rangle\Big|
\lesssim & 2^k \Big(\frac{2^n\delta_1}{\delta_0}\Big)^{N+1-j}, \qquad 1 \le j \le N, 
\\[3pt]
 \label{yN} 
(2B)^{-1} 2^{k-1}  \le \big|\big\langle \mathbf y_\mu(\tau,\xi), \mathbf r^{(N+1)}&
  \big\rangle\big|\le B 2^{k+1}.
 \end{align}
\end{lem}
 
 \begin{proof} We write $\mathbf r=(\mathbf r_1,\dots,\mathbf r_{N+1})$. Note $\mathbf r^{(j)}_m(s)=s^{m-j}/(m-j)!$ for $m\ge j$. By \eqref{formulaw} we have
 \begin{align*}
y^{m-1}_\mu \mathbf r^{(j)}_m(2^n\delta_1\nu-\delta_0\mu)
&=\sum_{\ell=m-1}^{N}\mathfrak g^{\ell}_\mu\frac{(2^n\delta_1\nu-\delta_0\mu)^{m-j}}{(\ell+1-m)!(m-j)!}\,\omega_\mu^{\ell+1-m} 
\end{align*}
for $m\ge j$.  Since $\mathbf r^{(j)}_m(s)=0$ for $j>m$, taking sum  over $m$ gives  
\begin{equation*}
\big\langle \mathbf y_\mu, \mathbf r^{(j)}(2^n\delta_1\nu-\delta_0\mu)\big\rangle
=\sum_{\ell=j-1}^{N}\mathfrak g^{\ell}_\mu\frac{( 2^n\delta_1\nu-\delta_0\mu+ \omega_\mu)^{\ell+1-j}}{(\ell+1-j)!}. 
\end{equation*}

From \eqref{s} note $2^n\delta_1\nu-\delta_0\mu +\omega_\mu=2^n\delta_1\nu-\sigma(\xi)+\mathcal E_1$. Thus, \eqref{si-2n}, \eqref{ejb} with $j=1$, and \eqref{ord2} with $j=1$
 show  $|2^n\delta_1\nu-\delta_0\mu+\omega_\mu |\lesssim 2^n\delta_1$. Using  \eqref{gbound}, 
  we obtain  
  \[ 
 \big|\big\langle \mathbf y_\mu(\tau,\xi), \mathbf r^{(j)}(2^n\delta_1\nu-\delta_0\mu)\big\rangle\big|
\lesssim 2^k (2^{n}\delta_1)^{N+1-j}, \qquad 1 \le j \le N. \]
By homogeneity   it follows that $\langle \eta,  \mathbf r^{(j)}(\delta_0 s)\rangle=\delta_0^{N+1-j}\langle \mathcal D_{\delta_0} \eta, \mathbf r^{(j)}( s)\rangle$ 
  for $\eta\in \mathbb R^{N+1}$. 
  Therefore, we get \eqref{yy0}.    For   \eqref{yN}   note $\mathbf r^{(N+1)}=(0,\dots, 0, 1)$, so $\langle \mathbf y_\mu, \mathbf r^{(N+1)}\rangle=y^{N}_\mu $  and   \eqref{yN} follows by  \eqref{YN}.
\end{proof}

 Let $\mathrm V= \sspan\{\gamma'(\delta_0\mu),\dots,\gamma^{(N)}(\delta_0\mu)\}$ and 
 $\{v_{N+1},\dots, v_d\}$ be an orthonormal basis of  $\mathrm V^\perp$.   
Since $\gamma$ satisfies  $\vcon{N}B$, for each  $\xi \in \mathbb R^d$ we can write 
\Be 
\label{coord} \xi= \overline \xi+\sum_{N+1\le j\le d}y_j(\xi) v_j,
\Ee
where $\overline \xi\in \mathrm V$
and $y_j(\xi) \in \mathbb R$, $N+1\le j\le d$. We  define a linear map  $\mathrm Y_\mu^{\delta_0}$   by    
\[  \mathrm Y_\mu^{\delta_0}(\tau,\xi)= \big( 2^{-k} \mathcal D_{\delta_0} \mathbf y_\mu(\tau, \xi),\, y_{N+1}(\xi),\dots,\,y_d(\xi)\big).
\] 
Then, by \eqref{yy0} and \eqref{yN} we see
\Be
\label{supp-supp}
\mathrm Y_\mu^{\delta_0} ( \supp_{\tau,\xi} \mathfrak a_{\nu}^{\mu,n} )
\subset \mathbf s\Big(\frac{2^n\delta_1}{\delta_0}\nu-\mu,   C\frac{2^n\delta_1}{\delta_0},  2^2B; \mathbf r_\circ^{N+1}\Big) \times \mathbb R^{d-N}
\Ee
for some $C>1$.  Thus, we have  the inequality  \eqref{dgain0} for $\delta=C{2^n\delta_1}/{\delta_0}$,  the collection of 
slabs $\mathbf s(2^n\delta_1\nu/\delta_0-\mu,  C2^{n}\delta_1/\delta_0, CB; \mathbf r_\circ^{N+1}), $ $\nu\in \mathfrak J_{n}^\mu$.  Therefore, via cylindrical extension  in $y_{N+1},\dots,y_d$ (see \eqref{cylinder}) and  the change of  variables $(\tau,\xi)\to \mathrm Y_\mu^{\delta_0}(\tau,\xi)$  we obtain \eqref{decoupn} since the decoupling inequality is not affected  by affine change of variables
in the Fourier side.

Combining \eqref{apnorm} and \eqref{decoupn}, we obtain 
\begin{equation*}
 \big(\sum_\mu \| \mathcal T[\mathfrak a^\mu]f\|_p^p\,\big)^{1/p}
	\le 
		\sum_{n\ge 0}\, \mathbf E_n
\end{equation*}
for $p \ge 4N-2$, 
where
\begin{align*}
\mathbf E_n= C_\epsilon
 \big({2^n\delta_1}/\delta_0\big)^{\frac {2N} {p}-1-\epsilon}  \big(\sum_{\mu} \sum_{\nu\in \mathfrak J_n^\mu} \| \mathcal T[\mathfrak a^{\mu, n}_\nu]f\|_p^p\,\big)^{1/p}.
\end{align*}

Since the intervals $I(\delta_0\mu,\delta_0)$  overlap, there are
at most three nonzero  $\mathfrak a^{\mu, 0}_{\nu}$ for each $\nu$.  We take  $\mathfrak a_\nu=\mathfrak a^{\mu, 0}_\nu$  which maximizes 
$\| \mathcal T[\mathfrak a^{\mu, 0}_\nu]f\|_p$.  Then,   it is clear that 
$\mathbf E_0\le  3^{1/p} C_\epsilon
 (\delta_1/\delta_0)^{\frac {2N} {p}-1-\epsilon} (\sum_{\nu} \| \mathcal T[\mathfrak a_\nu ]f\|_p^p\,)^{1/p}.$ By Lemma \ref{itsym},  $C^{-1}\mathfrak a_\nu\in \mathfrak A_k(\delta_1\nu, \delta_1)$ for a constant $C$. Thus, the proof of Proposition \ref{iterative}  is now reduced to showing 
\begin{align}\label{S*}
\sum_{n\ge1}\, \mathbf E_n \lesssim_B\! \delta_0^{-\frac{2N}p+1+\epsilon} 2^{-\frac 2pk +2\epsilon k} \|f\|_p,   \quad   p\ge 4N-2.
\end{align}

\subsubsection*{Estimates for $\mathbf E_n$ when $n\ge 1$} 
To show \eqref{S*} we decompose $\mathfrak a^{\mu, n}_\nu$  so that  \eqref{sumss} or \eqref{sumss2} (see Lemma \ref{lower-lower} below) holds on the $s,\xi$-supports of the resulting symbols. 
 If  \eqref{sumss} holds, we use the assumption  after rescaling, whereas we handle  the other case using estimates for the kernels of the operators.

 Let 
\Be
\label{G0}
\bar{\mathfrak G}_N^{\mu}(s,\xi)
=\sum_{1\le j\le N-2}\big(2^{-k}\mathfrak g^{j}_\mu\big)^{\frac{2N!}{N-j}}+\big(s-\sigma(\xi)\big)^{2N!}.
\Ee
The right hand side is independent of $\tau$ since so are $\mathfrak g^{j}_\mu$, $1\le j\le N-2$.

Let  $C_0=2^{2d}B$. We set  
\Be
\label{sym-def}
\mathfrak a_{\nu,1}^{\mu,n} =\mathfrak a_{\nu}^{\mu,n}  \,
\beta_0 \Big({ (2^{-k}\mathfrak g_\mu^0)^{2(N-1)!}}/
({C_0^{2N!} \bar{\mathfrak G}_N^{\mu}})\Big),  \quad n\ge 1,
\Ee
and  $\mathfrak a_{\nu,2}^{\mu,n}=\mathfrak a_{\nu}^{\mu,n} -\mathfrak a_{\nu,1}^{\mu,n},$ so  $ 
\mathfrak a_{\nu}^{\mu,n} =\mathfrak a_{\nu,1}^{\mu,n}+\mathfrak a_{\nu,2}^{\mu,n}.
$
Similarly as before,  we have the following, which we prove  in Section \ref{pf-lem2}. 
\begin{lem}\label{itsym2} 
There exists a constant $C=C(B)$ such that $C^{-1}\mathfrak a_{\nu,1}^{\mu,n}$, and $C^{-1} \mathfrak a_{\nu,2}^{\mu,n}$
 are contained in $\mathfrak A_k(2^n\delta_1\nu, 2^n\delta_1)$ for $n\ge 1$. 
\end{lem}

The estimate \eqref{S*}  follows if we show
\begin{align}\label{T1pp}
\big(
\sum_{\mu} \sum_{\nu \in \mathfrak J_n^\mu} \|\mathcal T[\mathfrak a_{\nu,1}^{\mu,n}]f\|_p^p \big)^{1/p}
\le C_\epsilon
2^{-\frac 2pk+\epsilon k}
(2^n\delta_1)^{-\frac{2N}p+1+\epsilon}
\|f\|_p,  \,  \quad p\ge 4N-6, \hspace{-4pt}
\end{align}
for any $\epsilon>0$, and 
\begin{align}
\label{Tpp}
\big( \sum_\mu\sum_{\nu\in \mathfrak J_{n}^{\mu}}\| \mathcal T[\mathfrak a_{\nu,2}^{\mu,n}]f\|_p^p \big)^{1/p}
\lesssim_B\!  2^{-\frac {(N+2)k} {2N}}(2^n\delta_1)^{-\frac N2}\|f\|_p, \quad  2\le p\le \infty  
\end{align} 
when $n\ge 1$.  Thanks to   \eqref{nnn}, those estimates  give 
\[
\sum_{n\ge1} \mathbf E_n\,\le C_\epsilon   \delta_0^{-\frac{2N}p+1+\epsilon}  \!\! \sum_{1\le n\le \log_2 (C\delta_0/\delta_1)}  
\!\!\big(2^{-\frac {2}pk + \epsilon  k}
+    2^{-\frac {(N+2)k} {2N}}(2^n\delta_1)^{\frac{2N}p-\frac {N+2}2-\epsilon}\big) \|f\|_p  \]
for $p\ge 4N-2$.  Note  $ \log_2 (\delta_0/\delta_1)\le Ck $ from \eqref{orderB}.  So,  \eqref{S*}  follows 
since $4N-2>  4N/(N+2)$ and $\delta_1\ge 2^{-k/N}$. 

\medskip

In order to prove the estimates \eqref{T1pp} and \eqref{Tpp},  we start  with the next lemma.  

\begin{lem}
\label{lower-lower} Let $n\ge 1$.
For a constant $C=C(B)>0$,  we have the following\,$:$ 
\begin{align}
\label{sumss}
\sum_{1\le j\le N-1} (2^n\delta_1)^{-(N-j)} |\langle\gamma^{(j)}(s),\xi\rangle| \ge C 2^k,  \qquad (s,\xi)\in  \supp_{s,\xi} \mathfrak a_{\nu,1}^{\mu,n}, 
\\
\label{sumss2}
(2^n\delta_1)^{-N} |\tau+\langle\gamma(s),\xi\rangle|\ge C 2^{k}, \qquad (s,\xi)\in   \supp_{s,\xi} \mathfrak a_{\nu,2}^{\mu,n}.
\end{align} 
\end{lem}

\begin{proof} We first prove \eqref{sumss}. Since $\mathfrak G_N^\mu\ge2^{-2N!-1}(2^n\delta_1)^{2N!}$ on  $\supp_{s,\xi} \mathfrak a_{\nu}^{\mu,n}$,  
one of the following holds on $\supp \mathfrak a_{\nu,1}^{\mu,n}$: 
\begin{align} 
\label{small}
|s-\sigma(\xi)|&\ge (2^{3}C_0B)^{-1}2^n\delta_1, 
\\
\label{lowergg}
2^{-k}|\mathfrak g^{j}_\mu(\tau,\xi)| &\ge  (2^2C_0)^{-(N-j)} (2^n\delta_1)^{N-j} 
\end{align}
for some $1\le j \le N-2$, where $C_0=2^{2d}B$ (see \eqref{sym-def}). 
If \eqref{small} holds, by \eqref{lowN} and \eqref{zero} it follows that   $
(2^n\delta_1)^{-1}|\langle\gamma^{(N-1)}(s),\xi\rangle| \gtrsim 2^{k}.$ 
Thus, to show \eqref{sumss} we may  assume  \eqref{small} fails, i.e.,  \eqref{lowergg} holds for  some $1\le j \le N-2$. So,  
 there is an  integer $\ell\in [0, N-2]$ such  that   \eqref{lowergg} fails for $\ell+1\le j \le N-2$,  whereas  \eqref{lowergg} holds for $j=\ell$. By \eqref{uj} and 
 \eqref{ejb}, we  have 
\Be
\label{haaa} |\langle G^{(\ell)}(\sigma(\xi)),(\tau,\xi) \rangle|
\ge |\mathfrak g_\mu^\ell| -\!\!\sum_{j=\ell+1}^N|\mathfrak g_\mu^j|
\frac{(B^6 2^{14}\delta_0^2)^{j-\ell}}{(j-\ell)!} -2B^3 (B^2 2^7\delta_0)^{N+1-\ell}|\xi|.
\Ee 

Thus, by  \eqref{ord2}, $|\langle G^{(\ell)}(\sigma(\xi)),(\tau,\xi) \rangle| \ge  
(2^3C_0)^{-(N-\ell)} 2^k (2^n\delta_1)^{N-\ell}$. Also,  \eqref{uj} and our choice of $\ell$ give
$|\langle G^{(j)}(\sigma(\xi)),(\tau,\xi) \rangle| \le 
(2C_0)^{-(N-j)} 2^k (2^n\delta_1)^{N-j}$ for $\ell+1 \le j \le N-2$. Combining this with  $|s-\sigma(\xi)|< (2^3 C_0B)^{-1} 2^n\delta_1$
and expanding $G^{(\ell)}$ in Taylor series at  $\sigma(\xi)$, 
we see  $ |\langle G^{(\ell)}(s), (\tau,\xi) \rangle|\ge  C 2^k (2^n\delta_1)^{N-\ell}$ for some $C=C(B)>0$. 
  This proves  \eqref{sumss}.

We now show  \eqref{sumss2}, which is easier. On $\supp \mathfrak a_{\nu,2}^{\mu,n}$, 
$2^{-k}|\mathfrak g_\mu^0|\ge 2^{-N-1}(2^n\delta_1)^N$ and $2^{-k}|\mathfrak g_\mu^j| \le 2C_0^{-(N-j)}(2^n\delta_1)^{N-j}$ for $j=1, \dots, N-2$. 
Using \eqref{haaa} with $\ell=0$, by \eqref{ord2} and \eqref{orderB}  
we get   $ (2^n\delta_1)^{-N} |\tau+\langle\gamma(\sigma(\xi)),\xi\rangle| \ge 2^{-N-2}2^k$. 
 We also note that $|s-\sigma(\xi)| \le 2C_0^{-1} 2^n\delta_1$ and  $|\langle G^{(j)}(\sigma(\xi)),(\tau,\xi) \rangle| \le C_0^{-1}2^k (2^n\delta_1)^{N-j}$ 
for $1 \le j \le N-2$ on $\supp \mathfrak a_{\nu,2}^{\mu,n}$. 
Since $|\langle G^{(N)}(s),(\tau,\xi) \rangle| \le B2^{k+1}$, using  Taylor series  expansion at $\sigma(\xi)$  as above,   we see  \eqref{sumss2}  holds true for some $C=C(B)>0$.  
\end{proof}

Additionally, we make use of disjointness of  $\supp_\xi\mathfrak a_{\nu}^{\mu,n}$ by  combining Lemma \ref{disjoint} and the next one. 

\begin{lem}
\label{L-supp}  There is a  positive constant $C=C(B)$ such that 
\Be  |( \widetilde{\mathcal L}_s^\delta)^{-1}\xi|\le Cb2^{k}
\label{haa}\Ee
whenever  $\xi\in {\Lambda}'_k(s,\delta,b)$ $($see \eqref{wl}$)$. If $\xi\in  \Gamma_k$ and \eqref{haa} holds with $C=1$, then  $\xi\in  {\Lambda}'_k(s,\delta,  C_1b)$ for some $C_1=C_1(B)>0$. 
\end{lem}

\begin{proof} Let $\eta\in \mathbb R^d$ and $\{ v_{N}, \dots, v_d\}$ be an orthonormal  basis of $(\sspan\{  \gamma^{(j)}(s): 1\le j\le N-1\})^\perp$.  We write 
$ \eta =\sum_{j=1}^{N-1}  \mathbf c_j \gamma^{(j)}(s)  + \sum_{j=N}^d \mathbf c_j v_j.$   Since $\vcon N B$ holds for $\gamma$, $|\eta|\sim |(\mathbf c_1, \cdots,\mathbf c_d )|$.  
Let $\xi\in {\Lambda}'_k(s,\delta,b) $.  Then,    \eqref{overL} gives
\[
\langle \eta, ( \widetilde{\mathcal L}_s^\delta)^{-1}\xi \rangle= \langle  ( \widetilde{\mathcal L}_s^\delta)^{-\intercal} \eta, \xi \rangle=
\sum_{j=1}^{N-1}  \delta^{j-N} \mathbf c_j \langle\gamma^{(j)}(s), \xi\rangle   + \sum_{j=N}^d \mathbf c_j \langle  v_j,\xi\rangle.  
\]
Thus, by $\eqref{wl}$  we get $|\langle \eta, ( \widetilde{\mathcal L}_s^\delta)^{-1}\xi \rangle|\le Cb|\eta| 2^k $, which shows \eqref{haa}. 

By \eqref{overL}, $\langle \gamma^{(j)}(s),\xi \rangle=  \delta^{N-j} \langle \gamma^{(j)}(s),  (\widetilde{\mathcal L}_s^\delta)^{-1}\xi \rangle $ for $1\le j\le N-1$. 
Therefore,  \eqref{haa} with $C=1$ gives  $|\langle \gamma^{(j)}(s),\xi \rangle| \le C_1  b  \delta^{N-j} 2^k$ for a constant $C_1>0$ when $1\le j\le N-1$. 
This proves the second statement. 
\end{proof}

Now,  we are ready to prove the estimates \eqref{T1pp} and  \eqref{Tpp}.  We first show \eqref{T1pp}.

\begin{proof}[Proof of  \eqref{T1pp}] 
By Lemma \ref{itsym2},  $C^{-1}\mathfrak a_{\nu,1}^{\mu,n}\in \mathfrak A_k(2^n\delta_1\nu, 2^n\delta_1)$  for some $C>0$, and  
\eqref{sumss} holds on $\supp_{s,\xi} \mathfrak a_{\nu,1}^{\mu,n}$. Thus,  taking $\delta=2^n\delta_1$ and $\sz= 2^n\delta_1\nu$,  
we may use Lemma \ref{lem:res} for $\tilde\chi \mathcal T [\fa_{\nu,1}^{\mu,n}] f$ to get 
\begin{align*}
\big\| \tilde\chi \mathcal T [\fa_{\nu,1}^{\mu,n}] f \big\|_{L^p(\mathbb R^{d+1})}
\le C\sum_{1\le l\le C}  \delta
\big\|  \mathcal A_t 
[\gamma_{\sz}^{\delta}, \,\raisebox{-.2ex}{$ a_l$} \,]
\tilde f_l \big\|_{L^p(\mathbb R^{d+1})},
\end{align*}
 where  $ \|\tilde f_l\|_p= \|f\|_p$, $a_{l}$ are of type $(j, N-1, B')$  relative  to $\gamma_{\sz}^\delta$ for some $B'>0$, and $2^j\sim (2^n\delta_1)^N 2^k$. 
As seen before, $ \gamma=\gamma_{\sz}^{\delta}$ satisfies  $\vcon{N}{3B}$ and \eqref{curveB} with $B$ replaced by $3B$  for $\delta\le \delta_*$.
  So,  $ \gamma=\gamma_{\sz}^{\delta}$ satisfies $\vcon{N-1}{B'}$  for a large $B'$.  

Therefore, we may apply the assumption (Theorem \ref{lclb} with $L=N-1$) to $ \mathcal A_t [{\gamma_{\sz}^{\delta}}, a_l] $,  
which gives 
$
\| \mathcal A_t [{\gamma_{\sz}^{\delta}}, a_l] f\|_{p
}\le
C_\epsilon \big(2^{k} (2^{n}\delta_1)^N \big)^{-\frac {2}p + \epsilon } \|f\|_{p
} $
for a constant   $C_\epsilon=C_\epsilon(B')$. 
Consequently, we obtain 
\[\|  \tilde\chi \mathcal T[\mathfrak a_{\nu,1}^{\mu,n}]f\|_p  \le C_\epsilon  2^{-\frac 2pk+\epsilon k}
(2^n\delta_1)^{1-\frac{2N}p+\epsilon}\|f\|_p\] 
 for $p\ge 4(N-1)-2$. Besides, since $C^{-1}\mathfrak a_{\nu,1}^{\mu,n}\in \mathfrak A_k(2^n\delta_1\nu, 2^n\delta_1)$,  by \eqref{kernelrmk} we  have 
$\| (1-\tilde\chi)\mathcal T[\mathfrak a_{\nu,1}^{\mu,n}]f\|_{L^p(\mathbb R^{d+1})} \lesssim_ B\!2^{-k}(2^n\delta_1)^{1-N} \|f\|_{L^p(\mathbb R^d)}$ for $p>1$. 
Note $2^n\delta_1\gtrsim 2^{-k/N}$. Combining  those two estimates yields 
\Be 
\label{scaled}\| \mathcal T[\mathfrak a_{\nu,1}^{\mu,n}]f\|_p  \le C_\epsilon 
2^{-\frac 2pk+\epsilon k}
(2^n\delta_1)^{1-\frac{2N}p+\epsilon}
\|f\|_p.
\Ee

To  exploit disjointness of $\supp_\xi\mathfrak a_{\nu,1}^{\mu,n}$, we define 
a multiplier operator  by
\[ \mathcal F(P_s^\delta f)(\xi)= \beta_0\big(| (\widetilde {\mathcal L}_s^\delta)^{-1}\xi|/(C_02^k) \big)\widehat f(\xi)\]
for a constant $C_0>0$. Since $\supp_\xi \mathfrak a_{\nu,1}^{\mu,n}\subset {\Lambda}'_k (2^n\delta_1\nu, $ $2^n\delta_1,  2^5B)$, 
by Lemma \ref{L-supp} we may choose   $C_0$ large enough so that $\beta_0\big(| (\widetilde {\mathcal L}_{2^n\delta_1\nu}^{2^n\delta_1})^{-1}\cdot |/(C_02^k) \big)=1 $  
on $\supp_\xi \mathfrak a_{\nu,1}^{\mu,n}$.  
Thus,  $\mathcal T[\mathfrak a_{\nu,1}^{\mu,n}]f=  \mathcal T[\mathfrak a_{\nu,1}^{\mu,n}]   P_{2^n\delta_1\nu}^{2^n\delta_1} f $. 
Combining this and \eqref{scaled},  we obtain 
\[ \big(\sum_{\mu}\sum_{ \nu \in \mathfrak J_n^\mu} \| \mathcal T[\mathfrak a_{\nu,1}^{\mu,n}]f\|_p^p \,\big)^{1/p}
\le  C_\epsilon  2^{-\frac 2pk+\epsilon k}
(2^n\delta_1)^{1-\frac{2N}p+\epsilon}  \big( \sum_{\mu}\sum_{ \nu \in \mathfrak J_n^\mu} \|P_{2^n\delta_1\nu}^{2^n\delta_1} f\|_p^p\, \big)^{1/p} 
\] 
for a constant   $C_\epsilon=C_\epsilon(B)$ if $p\ge 4N-6$.  Therefore,  \eqref{T1pp} follows if we show 
\Be
\label{lplp}
 \big(
\sum_{\mu}\sum_{ \nu \in \mathfrak J_n^\mu} \|P_{2^n\delta_1\nu}^{2^n\delta_1} f\|_p^p\,\big)^{1/p}\lesssim_B\! \|f\|_p, \qquad  2\le p\le\infty. 
\Ee

By interpolation it suffices to obtain  \eqref{lplp} for  $p=2, \infty$.  The case 
$p=\infty$ is trivial since $\|P_{2^n\delta_1\nu}^{2^n\delta_1} f\|_\infty\lesssim \|f\|_\infty$.  
For $p=2$,  \eqref{lplp} follows by Plancherel's theorem since 
$\supp \beta_0\big(| (\widetilde {\mathcal L}_{2^n\delta_1\nu}^{2^n\delta_1})^{-1}\cdot |/(C_02^k) \big)\widehat f$, $\nu \in \mathfrak J_n^\mu$ are finitely overlapping.  Indeed, by Lemma \ref{L-supp}  we have $\supp \beta_0\big(| (\widetilde {\mathcal L}_{2^n\delta_1\nu}^{2^n\delta_1})^{-1}\cdot |/(C_02^k) \big)\widehat f\subset {\Lambda}'_k (2^n\delta_1\nu, 2^n\delta_1,   C_1B)$ 
for a constant $C_1$. It is clear from lemma \ref{disjoint} that ${\Lambda}'_k (2^n\delta_1\nu, 2^n\delta_1,   C_1 B)$, $\nu \in \mathfrak J_n^\mu$ overlap at most $C=C(B)$ times. 
  \end{proof}

The proof of \eqref{Tpp} is much easier since we have a favorable estimate for the kernel of  $\mathcal T[\mathfrak a_{\nu,2}^{\mu,n}]$
thanks to the lower bound \eqref{sumss2}. 

\begin{proof}[Proof of \eqref{Tpp}]
Let 
\[
\mathfrak  b(s,t,\tau,\xi)=
i^{-1}(\tau+\langle\gamma(s),\xi\rangle)^{-1}\partial_t \mathfrak a_{\nu,2}^{\mu,n} (s,t,\tau,\xi).
\]
Then, integration by parts in $t$ shows $m[\mathfrak a_{\nu,2}^{\mu,n}]=m[\mathfrak  b]$.   Note \eqref{sumss2} holds and $C^{-1}\mathfrak a_{\nu,2}^{\mu,n} \in \mathfrak A_k(2^n\delta_1\nu, 2^n\delta_1)$ for a constant $C\ge 1$. Thus,  $\mathfrak a:=C^{-1}2^k (2^n\delta_1)^N \mathfrak b$
satisfies, with  $\delta=2^n\delta_1$ and $\sz=2^n\delta_1\nu$,  \eqref{symsupp} and \eqref{symineq2}    for $0\le j\le 1$, $0\le l\le 2N-1$,  $|\alpha|\le d+N+2$.  Applying  
\eqref{ker-est}, we obtain 
$
\|\mathcal T[\mathfrak a_{\nu,2}^{\mu,n}]f\|_\infty \lesssim_B\! 2^{-k}(2^n\delta_1)^{1-N}\|f\|_\infty.
$ 
Since $\delta_1\ge 2^{-k/ N}$, this gives 
\begin{align}\label{inf-inf}
\|\mathcal T[\mathfrak a_{\nu,2}^{\mu,n}]f\|_\infty \lesssim_B \! 2^{-\frac {(N+2)k} {2N}}(2^n\delta_1)^{-\frac N2}\|f\|_\infty.
\end{align}

By interpolation it is sufficient to show  \eqref{Tpp} for $p=2$.
Note $\|b(\cdot,t,\tau,\xi)\|_\infty+ \|\partial_s b(\cdot,t,\tau,\xi)\|_1\lesssim 2^{-k} (2^n\delta_1)^{-N}.$
Thus,  \eqref{lowN} and  van der Corput's lemma in $s$  give 
$
|m[\mathfrak a_{\nu,2}^{\mu,n}](\tau,\xi)|
\lesssim 2^{-k(1+N)/N} (2^n\delta_1)^{-N}.
$
Since $\supp_\xi \mathfrak a_{\nu,2}^{\mu,n}$ $\subset {\Lambda}'_k ( 2^n\delta_1\nu, 2^n\delta_1,  2^5B)$,  as before, we have $\mathcal T[\mathfrak a^{\mu,n}_{2,\nu}]f=  \mathcal T[\mathfrak a^{\mu,n}_{2,\nu}]   P_{2^n\delta_1\nu}^{2^n\delta_1} f $ with $C_0>0$ large enough.  Thus, by Plancherel's theorem 
\begin{align*}
\| \mathcal T[\mathfrak a_{\nu,2}^{\mu,n}] f \|_{L^2}^2
\lesssim_B\! 2^{-\frac{2(1+N)}N k} (2^n\delta_1)^{-2N}\!\!\!\iint_{\{ \tau: |\mathfrak g_\mu^0(\tau,\xi)| \le 2^{k+1}(2^n\delta_1)^N \}}\!\!\!\!\! d\tau\,
|\mathcal F(P_{2^n\delta_1\nu}^{2^n\delta_1} f)(\xi)|^2\, d\xi.
\end{align*}
Combining this and \eqref{lplp} yields  \eqref{Tpp} for $p=2$.
\end{proof}

\subsection{Proof of Lemma \ref{itsym}}
\label{pf-lem}
To simplify notations,  we denote   \[ \delc=2^n\delta_1, \qquad \sza=2^n\delta_1\nu\] 
for the rest of this section.  To prove Lemma \ref{itsym}, we   verify   \eqref{symsupp} and \eqref{symineq2} with $\mathfrak a=\mathfrak a_{\nu}^{\mu,n}$, $\delta=\delc,$ and $\sz=\sza$. 
The first is easy. In fact, since $\mathfrak a^\mu \in \mathfrak A_k(\delta_0\mu,\delta_0)$ and $\supp_s  \mathfrak a_{\nu}^{\mu,n}\subset I(\sza, \delc)$,  we only  need to show
\begin{align} 
\label{gu-bd}
 |\langle G^{(j)}(\sza),(\tau,\xi)\rangle| \le B 2^{k+5} \delc^{N-j}, \quad j=0,\dots,N-1 
\end{align}
on $\supp_{\tau, \xi} \mathfrak a_{\nu}^{\mu,n}$.  Using  \eqref{uj} and \eqref{gbound} together with \eqref{ord2} and \eqref{ejb}, one can easily obtain 
\Be
\label{gu-bd2}
|\langle G^{(j)}(\sigma(\xi)),(\tau,\xi)\rangle| \le 2^{k+1} \delc^{N-j}, \quad  j=0,\dots,N-1
\Ee
on $\supp_{\tau, \xi} \mathfrak a_{\nu}^{\mu,n}$. Expanding $\langle G^{(j)}(s),(\tau,\xi)\rangle$ in Taylor's series at $\sigma(\xi)$ gives \eqref{gu-bd} since \eqref{si-2n} holds. 

\medskip

We now proceed to show \eqref{symineq2} with $\mathfrak a=\mathfrak a_{\nu}^{\mu,n}$, $\delta=\delc,$ and $\sz=\sza$.  Since  $\mathfrak a_{\nu}^{\mu,n}$ consists of three factors
$\mathfrak a^\mu$, $\beta_N(\delc^{-2N!}\mathfrak G_N^\mu)$, and $\zeta(\delc^{-1}s-\nu)$, by Leibniz's rule it is sufficient to consider  the derivatives of  each of them.
The bounds on  the derivatives $\zeta(\delc^{-1}s-\nu)$ are clear. So, it suffices to show \eqref{symineq2}  for 
\[\mathfrak a \  =  \  \mathfrak a^\mu, \ \beta_N(\delc^{-2N!}\mathfrak G_N^\mu)\] with   $\delta=\delc$ and $\sz=\sza$ 
whenever $(\tau,\xi ) \in \supp \mathfrak a_{\nu}^{\mu,n} (s,t,\mathcal L_{\sza}^{\delc}\cdot)$. 

We  handle $\mathfrak a^\mu$   first. That is to say, we show 
\Be 
\label{claim33}
\big| \partial_s^{j}\partial^{l}_t\partial^{\alpha}_{\tau,\xi}
 \big( \mathfrak a^\mu (s,t,\mathcal L_{\sza}^{\delc}(\tau,\xi) 
 )\big) \big|
\lesssim_B \delc^{-j} 
 |(\tau,\xi)|^{-|\alpha|}, \qquad (j,l,\alpha)\in \mathcal I_N, 
 \Ee
for $(\tau,\xi ) \in \supp \mathfrak a_{\nu}^{\mu,n} (s,t,\mathcal L_{\sza}^{\delc}\cdot)$.    Since $\mathfrak a^\mu \in \mathfrak A_k(\delta_0\mu, \delta_0)$ and $|\sza-\delta_0\mu| \le \delta_0$,
we have
\begin{align}\label{center}
|\partial_s^{j}\partial^{l}_t\partial_{\tau,\xi}^{\alpha} \big(
\mathfrak a^\mu (s,t, \mathcal L_{\sza}^{\delta_0} (\tau,\xi) )\big) |
\lesssim_B  \delta_0^{-j} |(\tau,\xi)|^{-|\alpha|}, \qquad (j,l,\alpha)\in \mathcal I_N . 
\end{align}   
One can  show this using  \eqref{LLt}.  We consider  $
\mathcal U:=(\mathcal L_{\sza}^{\delta_0})^{-1} \mathcal L_{\sza}^{\delc}.
$
By \eqref{hoho} we have $|\,\mathcal U^\intercal z|\lesssim_B\!  |z|$ because  $|\delta_0^{-1} 2^n \delta_1| \lesssim_B  1$. Thus,  \eqref{center} gives
\begin{align*}
|\partial_s^{j}\partial^{l}_t\partial_{\tau,\xi}^{\alpha} \big(
\mathfrak a^\mu (s,t, \mathcal L_{\sza}^{\delta_0} \mathcal U(\tau,\xi) )\big) |
\lesssim_B  \delta_0^{-j}|\mathcal U(\tau,\xi)|^{-|\alpha|}
\end{align*}
for $(\tau,\xi ) \in \supp \mathfrak a_{\nu}^{\mu,n} (s,t,\mathcal L_{\sza}^{\delc}\cdot)$.   

Let $(\tau,\xi)\in \supp \mathfrak a_{\nu}^{\mu,n} (s,t, \mathcal L_{\sza}^{\delc}\cdot)$.  Then, $\mathcal L_{\sza}^{\delta_0} \mathcal U(\tau,\xi)=\mathcal L_{\sza}^{\delc}(\tau,\xi)
\in \Lambda_k(\sza, \delc, B)$, so $|\widetilde{\mathcal L}_{\sza}^{\delc}\,\xi  |  
\sim  |(\tau,\xi)|$ by Lemma \ref{2ksupport}. 
This and  \eqref{Lnorm} 
give  
\[
 |(\tau,\xi)|\sim  |\widetilde{ \mathcal L}_{\sza}^{\delc}\,\xi|  
\le | \mathcal L_{\sza}^{\delc}(\tau,\xi)| 
\le |\,\mathcal U(\tau,\xi)| 
\]
for $(\tau,\xi)\in \supp \mathfrak a_{\nu}^{\mu,n} (s,t, \mathcal L_{\sza}^{\delc}\cdot)$.
So, we obtain \eqref{claim33} since $\delc\lesssim \delta_0$. 

\medskip

We  continue to show \eqref{symineq2} for $\mathfrak a=\beta_N(\delc^{-2N!}\mathfrak G_N^\mu)$.  
Note  
$(2^{n}\delta_1)^{-2N!}
	\mathfrak G_N^\mu$
is a sum of
$(\delc^{-1}(s-\sigma(\xi)))^{2N!}$ and  $(\delc^{-(N-j)}2^{-k}\mathfrak g^{j}_\mu)^{{2N!}/{(N-j)}}$, $0\le j\le N-2$.  
Since the exponents ${2N!}/{(N-j)}$ are even integers, for the desired bounds   on  $\partial^{\alpha}_{\tau,\xi} (\beta_N(\delc^{-2N!}\mathfrak G_N^\mu))$ 
it suffices to show   the same bounds  on  the derivatives of  
\[  \delc^{-1}(s-\sigma(\xi)), \qquad \delc^{-(N-j)}2^{-k} \mathfrak g^{j}_\mu, \quad 0\le j\le N-2 .\]

The  bound on $\partial_\xi^\alpha \delc^{-1}(s-\sigma)$ is a consequence of \eqref{C2} and  the following lemma.   
To simplify notations,  we  denote 
\[\Xi=   {\mathcal L}_{\sza}^{\delc}(\tau, \xi),  \qquad \tilde \Xi= \widetilde{\mathcal L}_{\sza}^{\delc}\xi.\]

\begin{lem}\label{derivcheck}
 If $\Xi \in \supp_{\tau,\xi} \mathfrak a_{\nu}^{\mu,n}$,  then we have
\begin{align}
| \delc^{-1} \partial^{\alpha}_{\xi}  (\sigma(\tilde \Xi )) |\lesssim_B\! |\xi|^{-|\alpha|}, \qquad 1\le |\alpha|\le 2d+2.
\label{claim31}
\end{align}
\end{lem}

\begin{proof}
By  \eqref{zero}, $\gamma^{(N-1)} ( \sigma (\tilde \Xi) ) \cdot\tilde \Xi=0$.
Differentiation gives 
\begin{align}\label{Lsigma}
\gamma^{(N)}(\sigma(\tilde \Xi))\cdot\tilde \Xi \,\,\nabla_\xi (  \sigma(\tilde \Xi) )
+(\widetilde{\mathcal L}_{{\sza}}^{\delc})^\intercal\gamma^{(N-1)} (\sigma(\tilde \Xi))=0.
\end{align}  
Denote  $s=\sigma(\tilde \Xi)$. By  \eqref{overL}, $(\widetilde{\mathcal L}_{{\sza}}^{\delc})^\intercal\gamma^{(N-1)}(s)=\delc (\widetilde{\mathcal L}_{{\sza}}^{\delc})^\intercal(\widetilde{\mathcal L}_{s}^{\delc})^{-\intercal}\gamma^{(N-1)}(s)$.   
Since $|\sza-s| \le \delc$, i.e., \eqref{gbound2},  by Lemma  \ref{similarmx} we have  $|(\widetilde{\mathcal L}_{{\sza}}^{\delc})^\intercal\gamma^{(N-1)}\big( \sigma(\tilde \Xi ) \big)|\lesssim_B\! \delc$.   Besides, $|\gamma^{(N)}(\sigma(\tilde \Xi)) \cdot\tilde \Xi | \gtrsim |\tilde \Xi| \sim 2^k$ (see \eqref{lowN}).  
Thus, \eqref{Lsigma} and \eqref{C2} give 
\begin{align*}
|\nabla_\xi (\sigma(\tilde \Xi))| \lesssim_B \delc|\xi|^{-1}, \end{align*}
which proves  \eqref{claim31} with $|\alpha|=1$. 

We show  the bounds on the derivatives of higher order by induction. Assume that \eqref{claim31} holds true for $|\alpha|\le L$. Let $\alpha'$ be a multi-index such that $|\alpha'|=L+1$. Then,  
differentiating \eqref{Lsigma} and using  the induction assumption, one can  easily see 
$\gamma^{(N)}(\sigma(\tilde \Xi ))\cdot\tilde \Xi\,\,\partial_\xi^{\alpha'}\! ( \sigma(\tilde \Xi))= O(\delc |\xi|^{-L})$, by which we get  \eqref{claim31} for $|\alpha|=L+1$. Since $\sigma\in C^{2d+2}$, one can continue this as far as $L\le 2d+1$. 
\end{proof}

The proof of Lemma \ref{itsym} is now completed if we show 
\begin{align}\label{decomp03}
\big| 2^{-k} \partial^{\alpha}_{\tau,\xi} \big(\mathfrak g_\mu^\ell ( \Xi
) \big) \big|
\lesssim_B\! \delc^{N-\ell}2^{-k|\alpha|}, \qquad |\alpha| \le d+N+2
\end{align}
for $0\le \ell \le N-2$  whenever  $\Xi \in \supp \mathfrak a_{\nu}^{\mu,n} (s,t,\cdot)$.
To this end, we  use the following.

\begin{lem}\label{Glem}
For $ j=0,\dots,N$,  we  set 
\[
A_j=\delc^{-(N-j)} 2^{-k}  \big\langle G^{(j)}(\sigma(\tilde \Xi)), \Xi \big\rangle. 
\]
If $(\tau,\xi) \in \supp \mathfrak a_{\nu}^{\mu,n}(s,t, \mathcal L_{\sza}^{\delc} \cdot)$, then  for  $ j=0,\dots,N$  we have
\begin{align}\label{Ajderiv}
|\partial_{\tau,\xi}^\alpha A_j| \lesssim_B\! |(\tau,\xi)|^{-|\alpha|},  \qquad  	1\le  |\alpha|\le 2d+2.
\end{align}
\end{lem}

\begin{proof}   When $j=N$, the estimate \eqref{Ajderiv} follows by Lemma \ref{derivcheck} and   \eqref{C2}.  So, we may assume $j\le N-1$.  Differentiating  $A_j$,  we have 
 \[ \nabla_{\tau,\xi}A_j= B_j+ D_j,\]
 where
\begin{align*} 
B_j=\delc^{-1}
\big(0, \nabla_\xi ( \sigma (\tilde \Xi)) \big) A_{j+1},
\qquad 
D_j=\delc^{-(N-j)} 2^{-k}  (\mathcal L_{{\sza}}^{\delc})^\intercal G^{(j)}(\sigma(\tilde \Xi)).
\end{align*}
Note $(\mathcal L_{{\sza}}^{\delc})^\intercal G^{(j)}(\sza)=\delc^{N-j}G^{(j)}(\sza)$ for $0\le j\le N-1$.  Since $|\sza-\sigma( \tilde\Xi) 
| \lesssim \delc$,  similarly  as before,  Lemma \ref{similarmx}  and \eqref{hoho}  
give 
\Be 
\label{haha} |(\mathcal L_{\sza}^{\delc})^\intercal G^{(j)}(\sigma(\tilde \Xi 
))|\lesssim_B\! \delc^{N-j},  \qquad 0\le j\le N-1. 
\Ee 
By Lemma \ref{derivcheck} and \eqref{gu-bd2}, $| B_j|\lesssim  |\xi|^{-1}$.   Thus, for 
$\Xi \in \Lambda_k(\sza, \delc, B)$,  we have 
\[
|\nabla_{\tau,\xi}A_j| \lesssim_B\! |\xi|^{-1}+2^{-k} \lesssim_B\! |(\tau,\xi)|^{-1}, \quad  j=0,\dots,N-1, 
\]
For the second  inequality we use \eqref{C2}. This gives \eqref{Ajderiv} when $|\alpha|=1$. 

To show \eqref{Ajderiv} for $2\le |\alpha|\le 2d+2$, we use backward  induction. By  \eqref{zero}  we note    
$A_{N-1}=0$, so \eqref{Ajderiv} trivially  holds  
when $j=N-1$. We now assume that 
\eqref{Ajderiv} holds true  if $j_0+1\le j\le N-1$ for some $j_0$.   Lemma \ref{derivcheck}, \eqref{C2},  and  the induction assumption show $\partial^{\alpha'}_{\tau, \xi}  B_{j_0}=O(|(\tau,\xi)|^{-1-|\alpha'|})$  for  $1\le |\alpha'|\le 2d+1$.
 Concerning $ D_{j_0}$,  observe that  $\partial^{\alpha'}_{\xi} ( G^{(j_0)}(\sigma(\tilde \Xi 
 )))$ is given by a sum  of the terms
\[
G^{(j)}(\sigma(\tilde \Xi 
))\prod_{n=1}^{j-j_0}\partial^{\alpha'_{n}}_{\xi}(\sigma(\tilde \Xi 
)), 
\]
where $j\ge j_0$ and $\alpha'_1+\cdots+\alpha'_{j-j_0}=\alpha'$.   
Hence,   Lemma \ref{derivcheck}, \eqref{haha}, and  \eqref{C2}  give  $\partial^{\alpha'}_{ \xi}  D_{j_0} =O(|(\tau,\xi)|^{-1-|\alpha'|})$  for  $1\le |\alpha'|\le 2d+1$.  
Therefore, combining the estimates for $B_{j_0}$ and $D_{j_0}$, we get $\partial^{\alpha'}_{\tau, \xi} \nabla_{\tau,\xi}A_{j_0}= O(|(\tau,\xi)|^{-1-|\alpha'|})$. This proves \eqref{Ajderiv}  for $j=j_0$. 
\end{proof}

Before proving  \eqref{decomp03}, we first note   
\Be\label{d-ej}
 |\partial_\xi^\alpha\big(\mathcal E_j(\tilde \Xi )\big)|\lesssim_B \delc^j|\xi|^{-|\alpha|}, \quad |\alpha|\le 2d+2
 \Ee
 for $j=1, \dots, N$. This can be shown  by  a routine  computation. Indeed,  differentiating \eqref{ej}, and using  Lemma \ref{derivcheck} and   \eqref{ord2},  
 one can easily see \eqref{d-ej} since $|\sigma(\tilde \Xi)-\delta_0\mu|\lesssim \delta_0$.  

To show \eqref{decomp03} for $0\le \ell \le N-2$, we  again use  backward induction. Observe that  \eqref{decomp03} holds  for $\ell=N,N-1$, and  assume  that \eqref{decomp03} holds for $j+1\le \ell\le N$ for some $j\le N-2$.  By \eqref{uj} we have  \[
2^{-k}\mathfrak g^{j}_\mu 
	= \delta_\ast^{N-j} A_j-\sum_{j+1\le \ell \le N} ( 2^{-k} \mathfrak g^{\ell}_\mu)
	(\mathcal E_1)^{\ell-j}/{(\ell-j)!}  +2^{-k}y^N_\mu \mathcal E_{N-j}.\] 
Thus, by Lemma \ref{Glem}  and \eqref{d-ej}, we get   \eqref{decomp03} with $\ell=j$.    This completes the proof of Lemma  \ref{itsym}.

\subsection{Proof of Lemma \ref{itsym2}} 
\label{pf-lem2}
Lemma \ref{itsym2} can be proved in a similar way as the previous subsection. So, we shall be brief. 

By Lemma \ref{itsym} we have $C^{-1} \mathfrak a_{\nu}^{\mu,n}\in \mathfrak A_k(\sza,\delc)$ for  a constant $C\ge1$, so it  suffices to show 
$C^{-1}\mathfrak a_{\nu,1}^{\mu,n}\in \mathfrak A_k(\sza,\delc)$ for some $C\ge 1$.  The support condition \eqref{symsupp}  is obvious, so
we need only to show  \eqref{symineq2} with $\mathfrak a=\mathfrak a_{\nu,1}^{\mu,n}$, $\delta=\delc$, and $\sz=\sza$. 
Moreover,  recalling \eqref{sym-def},  it is enough to  consider the additional factor only, i.e., to show  
\[
\Big| \partial_{\tau,\xi}^\alpha \Big( \beta_0 \bigg( \frac{ \big( \delc^{-N}2^{-k}\mathfrak g_\mu^0( \mathcal L^{\delc}_{\sza}(\tau,\xi))\big)^{2(N-1)!}}
{ C_0^{2N!}\delc^{-2N!} \bar{\mathfrak G}_N^{\mu}(s, \widetilde{\mathcal L}^{\delc}_{\sza} \xi )} \bigg)\Big) \Big|
\lesssim |(\tau,\xi)|^{-|\alpha|}
\]
for   $(\tau, \xi)\in \supp \mathfrak a_{\nu,1}^{\mu,n}(s,t,\mathcal L^{\delc}_{\sza}\cdot)$.
Since $\delc^{-2N!}\bar{\mathfrak G}_N^{\mu} \gtrsim 1$ on  $\supp_{s,\xi}\mathfrak  a_{\nu,1}^{\mu,n}$, 
one can obtain the estimate  in the same way as in the proof of Lemma \ref{itsym}.

\subsection{Sharpness of Theorem \ref{LS}}
Before closing this section, we show the optimality of the regularity  exponent $\alpha$ in Theorem  \ref{LS}.

\begin{prop}\label{prop:n}
Suppose \eqref{ls} holds for $\psi(0)\neq 0$. Then
$\alpha\le 2/p$.
\end{prop}

\begin{proof} 
We write $\gamma=(\gamma_1, \dots, \gamma_d)$. 
Via an affine change of variables, we may assume $\gamma_1(0)=0$ and $\gamma_1'(s) \neq0$  on an interval 
$J=[-\delta_0,\delta_0] $ for $0<\delta_0\ll 1$.  Since  $\psi(0)\neq 0$, we may also assume $\psi \ge 1$ on $J$.

We choose $\zeta_0 \in \mathcal S(\mathbb R)$ such that $\supp \widehat \zeta_0 \subset [-1,1]$  
and $ \zeta_0 \ge 1$ on $[-r_1,r_1]$
where $r_1=1+2\max \{ |\gamma(s)| : s \in J\}$. Denoting $\bar x = (x_1,\dots,x_{d-1})$ and $\bar \gamma(t)=(\gamma_1(t),\dots,\gamma_{d-1}(t))$,  we define  
\begin{align*} 
\bar{\mathcal A}_t h(x) =\int e^{it\lambda \gamma_d(s)} 
\zeta_0(x_d-t\gamma_d(s)) h(\bar x - t \bar \gamma(s)) 
\psi(s)\,ds. 
\end{align*} 

Let $\zeta \in \mathrm C_c^\infty((-2,2))$ be a positive function  such that $\zeta=1$ on $[-1,1]$. For  a positive constant $c\ll \delta_0$,  let $g_1(\bar x)=\sum_{\nu \in  \lambda^{-1}\mathbb Z \cap  [-c,c]} \zeta(\lambda|\bar  x+\bar\gamma(\nu)|)$. We consider 
\[
g(\bar x)=e^{-i\lambda \varphi(x_1)}g_1(\bar x), 
\]
where $\varphi(s)= \gamma_d \circ (- \gamma_1)^{-1}(s)$.   We  claim that,   if $c$ is small enough, 
\Be 
\label{a-lower}
|\bar{\mathcal A}_t g( x)| \gtrsim1, \qquad  (x,t)\in S_c, 
\Ee
where $S_c= \{  (x,t): |\bar x| \le  c \lambda^{-1}, \, |x_d| \le c, \,  |t-1| \le  c \lambda^{-1}\}$. 
To show this, note
\begin{align*} 
\bar{\mathcal A}_t g(x) =\int
 e^{i\lambda ( t\gamma_d(s)- \varphi (x_1-t\gamma_1(s)))} \zeta_0(x_d-t\gamma_d(s)) g_1(\bar x - t \bar \gamma(s))\,\psi(s)\,ds.
\end{align*}

Let $(x,t)\in S_c $. Then, $\supp g_1(\bar x - t \bar \gamma(\cdot)) \subset [-C_1c, C_1c]$ for some $C_1>0$. 
Since $\varphi(s)= \gamma_d \circ (- \gamma_1)^{-1}(s)$, by the mean value theorem  we see
$
|\varphi(x_1-t\gamma_1(s))  
-\gamma_d(s)|\le 2r_0c\lambda^{-1} 
$
where $r_0=10 r_1\max\{ |\partial_s \varphi (s)| : \, s \in (-\gamma_1)(J_\ast)\} $ and $J_\ast=[-(C_1+1)c, (C_1+1)c]$. 
Thus,  we have
\Be
\label{phase}
|t \gamma_d(s)- \varphi (x_1-t \gamma_1(s))| \le 3r_0c\lambda^{-1}. 
\Ee
Besides,  if $\lambda$ is sufficiently large, 
$
g_1(\bar x-t\bar \gamma(s))=\sum_{\nu \in \lambda^{-1}\mathbb Z \cap  [-c,c]} 
\zeta (\lambda|\bar x-(t-1) \bar\gamma(s)+\bar\gamma(\nu)-\bar \gamma(s)|)
\gtrsim 1
$
if $s \in [-c/2, c/2]$. Since  $ \supp g_1(\bar x - t \bar \gamma(\cdot))\subset J$ with $c$ small enough and $\zeta_0 (x_d-t\gamma_d(s))\ge 1$, we get
$\int \zeta_0 (x_d-t\gamma_d(s)) g_1(\bar x-t\bar \gamma(s)) \psi(s)\,ds \gtrsim 1.  
$
Therefore,  \eqref{a-lower} follows  by \eqref{phase}    if  $c$ is small enough, i.e., $c\ll 1/(3r_0)$.

We  set $f(x)=e^{-i\lambda x_d} \zeta_0(x_d) g(\bar x).$ Then, 
$ 
\chi(t) \mathcal A_t f(x)= e^{-i\lambda x_d}  \chi(t) 
\bar {\mathcal A}_t   g(x). 
$
By our choice of $\zeta_0$, $\supp \widehat f\subset \{\xi:   |\xi_d+\lambda|\le 1\}$, so  $\supp \mathcal F(\chi(t)\mathcal A_t f)\subset \{(\tau, \xi):   |\xi_d+\lambda|\le 1\}$.  This gives
\Be
\label{sobol}
\lambda^{\alpha}\|  \chi(t) \mathcal A_t f\|_{L^p(\mathbb R^{d+1})}
\lesssim \|  \chi(t) \mathcal A_t f \|_{L_\alpha^p(\mathbb R^{d+1})}. 
\Ee
Indeed, 
$
 \lambda^\alpha  \|  \chi(t) \mathcal A_t f  \|_{L^p(\mathbb R^{d+1})}
\lesssim \|  \chi(t) \mathcal A_t f \|_{L^p( \mathbb R_{t,\bar  x}; L_\alpha^p(\mathbb R_{x_d}))}$  by  Mihlin's multiplier theorem in $x_d$. 
Similarly, one also sees  $ \|F \|_{L^p( \mathbb R_{t,\bar  x}; L_\alpha^p(\mathbb R_{x_d}))}  \le C \|F \|_{L_\alpha^p(\mathbb R^{d+1})}  $
for $\alpha\ge 0$ and any $F$. Combining those inequalities gives \eqref{sobol}.   

From \eqref{a-lower}  we have $\|\chi(t) \mathcal A_t f\|_p= \|  \chi(t) \bar{\mathcal A}_t g \|_{p}  \ge C\lambda^{-d/p}.$  
Note that $\supp g$ is contained in a $O(\lambda^{-1})$-neighborhood of  $-\bar{\gamma}$, so it follows that  $\|f\|_p \lesssim \lambda^{-(d-2)/p}$.  Therefore, by \eqref{sobol}   the inequality  \eqref{ls} implies $\lambda^\alpha \lambda^{-d/p} \lesssim \lambda^{-(d-2)/p}$.  Taking $\lambda\to \infty$ gives $\alpha\le 2/p$. 
\end{proof}

\section{$L^p$ Sobolev regularity}\label{sec:sobolev}
In this section we prove  Theorem \ref{thm:improving}, whose proof  proceeds in a similar way as  that  of Theorem \ref{LS}. 
However, we provide some details to make it clear 
how the optimal bounds are achieved. 
Since there are no  $t$, $\tau$ variables for the symbols,   the proof is consequently simpler but some modifications are necessary. 

For a large $B\ge 1$,  we assume  
\Be
\label{curveB-1}
 \max_{0\le j\le 2d} 
  |\gamma^{(j)}(s)|\,\le\,  B,   \qquad s\in I. 
 \Ee
   Let $2 \le L \le d$.  For $\gamma$ satisfying  \ref{lindepN}  
 we say $ \bar a  \in \mathrm C^{d+1}(\mathbb R^{d+1})$   is  \emph{a symbol of type $(k,L, B)$ relative to $\gamma$} if 
 $ \supp \bar a \subset I\times \mathbb A_k$,   \ref{sumN} holds for $\gamma$ on $\supp \bar a$, and 
\begin{align} 
\label{a-sym}
|\partial_s^j \partial_\xi^\alpha  \bar a (s,\xi)  |\le B|\xi|^{-|\alpha|}
\end{align}
for $0 \le j \le 1$ and $|\alpha| \le d+1$.  As before, Theorem \ref{thm:improving} is a straightforward consequence of the following. We denote   $\mathcal A[\gamma, \bar a  ]=\mathcal A_1[\gamma, \bar a  ]$.

\begin{thm}\label{sob-induct}
Suppose $\gamma \in \mathrm C^{2d}(I)$ satisfies \eqref{curveB-1}  and $ \a $ is a symbol of type $(k,L,B)$ relative to $\gamma$ for some  $B\ge 1$. 
Then, if  $p> 2(L-1)$,  for a constant $C=C(B)$ 
\begin{align} 
\label{lp-regular}
\| \mathcal A[\gamma, \bar a  ]f \|_{L^p(\mathbb R^d)}
\le C  2^{-k/p} \|f\|_{L^p(\mathbb R^d)}.
\end{align}
\end{thm}

In order to prove  Theorem \ref{thm:improving}, we consider  $ \bar a_k (s,\xi):=\psi(s)\beta(2^{-k}|\xi|)$, where
$\beta \in \mathrm C_c^\infty((1/2,4))$.  By \eqref{nonv}  $ \bar a_k $ is a symbol of type $(k,d, B)$ relative  to $\gamma$ for some $B$, thus  Theorem \ref{sob-induct}  gives \eqref{lp-regular} for  $p> 2(d-1)$. The estimate \eqref{lp-regular}  for each dyadic pieces can be put together by the result  in \cite{PRS}. So,  we get \eqref{A1} with $\alpha=\alpha(p)$ when $p> 2(d-1)$ (e.g., see \cite{BGHS2}).

Interpolation with   $\|  \mathcal A[\gamma, \bar a _k ] f\|_{2} \lesssim 2^{-k/d} \|f\|_{2}$ which follows from \eqref{f-decay} gives  $\|  \mathcal A[\gamma,\bar a _k ] f\|_{p} \lesssim_B\! 2^{-\alpha k } \|f\|_{p}$ for $\alpha\le \alpha(p)$ with strict inequality when $p\in (2, 2(d-1)]$. 
Using those estimates, we can prove  Corollary \ref{type}. Indeed, if $\gamma$ is a curve of maximal type $\ell>d$,   a typical scaling argument gives    
$\|  \mathcal A[\gamma,\bar a _k ] f\|_{p} \lesssim_B\! 2^{-\min(\alpha(p), 1/\ell)k } \|f\|_{p}$  for $p\neq \ell$ when $ \ell\ge 2d-2$,  and 
for $p\in  [2,  2\ell/(2d-\ell))\cup ( 2d-2, \infty)$  when $d< \ell< 2d-2$. As above,  one can combine the estimates (\cite{PRS}) to get  \eqref{A1}.

\subsection{Proof of Theorem \ref{sob-induct}}    
The case $L=2$ is easy. Since  $ \a $ is a symbol of type $(k,2,B)$ relative to $\gamma$, van der Corput's lemma  and Plancherel's theorem give \eqref{lp-regular} for $p=2$.    
Interpolation with $L^\infty$ estimate  shows  \eqref{lp-regular} for $p\ge 2$. When $L\ge 3$,  we have the following, which 
immediately yields  Theorem \ref{sob-induct}.

\begin{prop}\label{newpropsob} Let $3 \le N \le d$.
Suppose Theorem \ref{sob-induct} holds for $L=N-1$.
Then Theorem \ref{sob-induct} holds true  with $L=N$.
\end{prop}
  
   To prove the proposition, we fix $N\in [3, d]$ and 
$\gamma$ satisfying $\vcon N B$,  and   
$\a$ of type $(k, N, B)$ relative to $\gamma$.

For $\szz$ and $\delta>0$ such that  $I(\szz, \delta)\subset I$,  let 
\[
\bar{\Lambda}_k(\szz, \delta, B) =\bigcap_{1\le j\le N-1} \big\{ \xi \in \mathbb A_k : | \langle \gamma^{(j)}(\szz), \xi \rangle| \le B2^{k+5} \delta^{N-j}\big \}.
\]
By  $\bar{\mathfrak A}_k(\szz,\delta)$  we denote  a collection  of  $ \bfa  \in \mathrm C^{d+1}(\mathbb R^{d+1}) $ such that    $\supp \bfa\subset I(\szz, \delta)\times \bar{\Lambda}_k(\szz, \delta, B)$ and
$
|\partial_s^{j} \partial^{\alpha}_{\xi} \bfa  (s,  \widetilde{\mathcal L}^{\delta}_{\szz}\xi)|
\le B\delta^{-j}2^{-k|\alpha|}
$
for $0 \le j \le 1$ and $|\alpha| \le d+1$.  

 The next lemma which plays the same role as Lemma \ref{lem:res}  can be shown through  routine adaptation of the proof of Lemma \ref{lem:res}.

\begin{lem}\label{resc2}
Let $  \bfa \in \bar{\mathfrak A}_k(\szz,\delta)$ and $j_\ast=\log (2^k \delta^N)$. 
Suppose    \eqref{lowD} holds on $\supp \bar{\mathfrak a}$. Then, there exist  constants $C$, $\tilde B\ge 1$, and  $\delta'>0$  depending on $B$, and  symbols 
$\bfa_{1}, \dots, \bfa_{l_\ast}$ of type $(j, N-1, \tilde B)$  relative to $\gamma_{\szz}^\delta$,   
such that 
\[ 
\| \mathcal A[\gamma,\bfa  ]f  \|_{L^p(\mathbb R^d)} \le C \delta \sum_{1\le l\le C}  
\big\| \mathcal A 
[\gamma_{\szz}^{\delta}, \raisebox{-.2ex}{$\bfa _l$}]
\tilde f_l \big\|_{L^p(\mathbb R^{d})}, 
\] 
$ \|\tilde f_l\|_p= \|f\|_p$,  and  $j\in [j_\ast-C,  j_\ast +C]$  as long as  $0<\delta<\delta'$. 
\end{lem}

The order of necessary regularity on $\gamma$ is reduced since $\bfa$ is independent of $\tau, t$.  
Actually,  we  may take $\tilde a (s,\xi)=\bfa(\delta s +\szz,  \delta^{-N}\widetilde{\mathcal L}_{\szz}^\delta \xi )$ while following the  \emph{Proof of Lemma \ref{lem:res}}   since validity of  \eqref{a-sym} is clear for $\a=\tilde a$.

 Using $\eta_N$ (see \eqref{a-N}), 
we break \[\mathcal A[\gamma,\bar a]=\mathcal A[\gamma, \a\eta_N]+  \mathcal A[\gamma, \a(1-\eta_N)].\] 
 Note that  $C^{-1} \a(1-\eta_N)$ is of type $(k, N-1, B')$ relative to $\gamma$ for some large constants $B'$ and $C$, so we may apply the assumption to $\mathcal A[\gamma, \a(1-\eta_N)]f$. 
Consequently,  we have the estimate \eqref{lp-regular} for $\bar a=\a(1-\eta_N)$ if $ p> 2N-4$. 

 To obtain the estimate for $\mathcal A[\gamma, \a\eta_N]$, as before, we may assume $\supp \a\eta_N\subset I(\sz, \delta_\ast)\times  \bar\Gamma_k$ for some $\sz$ and a small $\delta_\ast$. 
 Here, $\bar\Gamma_k$ is defined in the same way as $\Gamma_k$ by replacing  $a\eta_N$ by $\a\eta_N$.  
 Since \eqref{lowN} holds on $\supp (\a\eta_N)$, we may  work under the same 
\emph{Basic assumption} as in Section \ref{sec:prop2.3}. 
That is to say,  we have $\sigma$ on $\bar\Gamma_k$ satisfying  \eqref{zero} and  $\sigma(\xi) \in  I(\sz, \delta_\ast) $ for $\xi \in \bar\Gamma_k$. 
Furthermore,  $\sigma\in \mathrm C^{d+1}$  since $\gamma \in \mathrm C^{2d}(I)$, and  \eqref{decaying-sigma} holds for $\xi \in \bar\Gamma_k$ and   $|\alpha|\le d+1$.  Thus, \eqref{a-sym} remains valid 
for the symbols   given subsequently  by decomposing $\bar a$ 
with cutoff functions associated with $\sigma$, and  $\bar{\mathfrak G}_N^{\mu}.$

 Apparently, $C^{-1}\a\eta_N \in \bar{\mathfrak A}_k(\szz,\delta_*)$ for a constant  $C=C(B, \delta_\ast)$,  therefore 
 the proof of Proposition \ref{newpropsob} is completed if we show the following.

\begin{prop}\label{main2}
Let $3 \le N \le d$ and $\bfa \in \bar{\mathfrak A}_k (\szz,\delta_*)$
with $ \supp_\xi \bfa\subset \bar\Gamma_k$. 
Suppose Theorem \ref{sob-induct} holds for $L=N-1$. 
Then, if $p> 2(N-1)$, we have \eqref{lp-regular}.
\end{prop}

We prove  Proposition \ref{main2} using the next,  which 
corresponds to Proposition \ref{iterative}. In what follows,  
we denote $\mathcal A[\a]=\mathcal A[\gamma, \a]$.

\begin{prop}\label{iterative2} 
Let $\delta_0$ and  $\delta_1$ satisfy  \eqref{orderB}.  For $\mu$ such that  $\delta_0\mu \in   I(\sz, \delta_\ast)\cap \delta_0\mathbb Z$, let $ \bfa^\mu \in \bar{\mathfrak A}_k(\delta_0\mu, \delta_0)$ with $ \supp \bfa^\mu \subset I(\sz, \delta_\ast)\times \bar \Gamma_k$.  Suppose Theorem \ref{sob-induct} holds for $L=N-1$. Then, if $p\in (2N-2, \infty)$,  there are constants  $\epsilon_0>0$, $C_0=C_0(\epsilon_0, B)\ge2$, and symbols
 $ \bfa_\nu \in \bar{\mathfrak A}_k( \delta_1\nu,\delta_1)$  with 
$ \supp \bfa_\nu \subset I(\sz, \delta_\ast)\times\bar\Gamma_k$, $\nu \in \cup_\mu\mathfrak J_0^\mu$, 
such that
\begin{equation*} 
\big(\sum_\mu \| \mathcal A[\bfa ^\mu] f \|_{p}^p \,\big)^{\frac1p}
\le C_0
\big({\delta_1} / {\delta_0}\big)^{\frac Np-1+\epsilon_0}
\big(\sum_{\nu}
\| \mathcal A[\bfa _\nu]f \|_{p}^p \,\big)^{\frac1p}
+C_0 \delta_0^{-\frac Np+1} 2^{-\frac kp}\|f\|_p.
\end{equation*}
\end{prop}

Let  $\delta'$ be given as in Lemma \ref{resc2}, and let $\delta_\circ>0$ be a positive constant  such that 
\Be\label{del-zero2}
\delta_\circ \le \min \{ \delta', \,(2^{7d}B^6)^{-N} C_0^{-2N/\epsilon_0}\}. 
\Ee

\begin{proof}[Proof of Proposition \ref{main2}]
Set $\delta_0=\delta_\circ$,  and let $\delta_1,\dots, \delta_J$ be given by \eqref{defn-j}. 
Then, applying Proposition \ref{iterative2} iteratively up to $J$-th step  (see Section \ref{sec:prop2.5}), we have symbols $ \bfa_\nu \in \bar{\mathfrak A}_k( \delta_{\!J}\nu, \delta_{\!J})$, $\delta_{\!J}\nu\in  I(\szz, \delta_0)$, such that
\[
\big\| \mathcal A[\bfa ]f \big\|_{p}
\le  C_0^J\delta_{\!J}^{\frac {N}{p}-1+\epsilon_0}
	 \big(\sum_\nu \| \mathcal A[\bfa _{\nu}]f\|_p^p \big)^{1/p}
\!+ 2^{-\frac{k}p} \delta_0^{-\frac Np+1-\epsilon_0}\!   \sum_{0\le j\le J-1}C_0^{j+1}  \delta_j^{\epsilon_0}  \|f\|_p.
\] 
By \eqref{del-zero2} and \eqref{defn-j},  
$ \delta_j \le C_0^{-2((N+1)/N)^jN/{\epsilon_0} }$, $0 \le j \le J-1.$
So,  $\sum_{j=0}^{J-1} C_0^{j+1} \delta_j^{\epsilon_0} \le  C_1$ 
for a constant $C_1$,  and   $C_0^J\delta_J^{\epsilon_0} \le  C_1$.
Thus,  the matter is now reduced to showing 
\[ 
 \big(\sum_{\nu} \| \mathcal A[\bfa_\nu] f\|_{L^p(\mathbb R^d)}^p \,\big)^{1/p}
\lesssim_B\!  2^{-\frac kN}\|f\|_{L^p(\mathbb R^d)}, \qquad  2\le p\le \infty,
\] 
which corresponds to  \eqref{est03}.  The case $p=\infty$  follows from the estimate  
$\| \mathcal A[\bfa  ]f\|_{L^\infty}$ $\le C\delta \|f\|_{L^\infty}$ when $ \bfa  \in \bar{\mathfrak A}_k(\szz,\delta)$ for some $\szz, \delta$ (cf. \eqref{ker-est}). 
One can obtain  this  in the same manner as in the proof of Lemma \ref{kernel}. 
The  case $p=2$ can be handled similarly  as before, using Plancherel's theorem and van der Corput's lemma combined with   Lemma \ref{disjoint} and \eqref{lowN}.   
\end{proof}

The proof of Proposition \ref{iterative2} is similar to that of Proposition \ref{iterative}. 
Instead of \eqref{dgain0} we use the estimate \eqref{dgain}, in which  the exponent 
is adjusted to  the sharp Sobolev regularity estimate.  However, a similar approach breaks down if one tries  to obtain  the local smoothing estimate  
\eqref{ls} with the optimal regularity $\alpha=2/p$. To do so,  we need the inequality \eqref{c-decoupling2} for $4N-2<p \le N(N+1)$. However, there is no such estimate available when $N=2$.

\subsection{Proof of Proposition \ref{iterative2}}  
\label{pf:thm1}
Let $ \bfa^\mu \in \bar{\mathfrak A}_k(\delta_0\mu,\delta_0)$. For $\nu \in \mathfrak J_n^\mu$,   set  \[
 \bfa_{\nu}^{\mu,n} =  \bfa^\mu 
\times
\begin{cases}
	\beta_0 \big( 
	\delta_1^{-2N!} \, \bar{\mathfrak G}_N^{\mu} 
	\big) \, \zeta(\delta_1^{-1}s-\nu),   &   \ n=0,
	\\[2pt]
	\beta_N \big( (2^{n}\delta_1)^{-2N!}
	\, \bar{\mathfrak G}_N^{\mu} 
	\big) \,	\zeta(2^{-n}\delta_1^{-1}s-\nu),  &  \ n\ge 1,
	\end{cases}
\] 
 (see  \eqref{G0}).   Let $\bar {\mathbf y}_\mu= (y_\mu^1,\dots,y^{N}_\mu)$, and let $\bar{\mathcal D}_\delta$ denote the $N\times N$ matrix $(\delta^{1-N} \bar e_1, $ $\delta^{2-N} \bar e_2, \dots,$ $\delta^0 \bar e_{N})$ where $\bar e_j$ is the $j$-th standard  unit vector in $\mathbb R^{N}$.
 Recalling \eqref{coord}, we consider a linear map 
 \[  \bar{\mathrm Y}_\mu^{\delta_0}(\xi)= \big( 2^{-k} \bar{\mathcal D}_{\delta_0} \bar{\mathbf y}_\mu,\, y_{N+1},\dots,\,y_d\big). 
\] 

Let $\mathbf r$ denote the curve $\mathbf r_\circ^N$.  Note that   \eqref{gbound2} and \eqref{gbound}  hold on $\supp \bfa_{\nu}^{\mu,n}$. Similarly as  in  \emph{Proof of Lemma \ref{supp-a}},  we see  $ | \langle \bar{\mathbf y}_\mu , \mathbf r^{(j)}((2^n\delta_1/\delta_0)\nu-\mu) \rangle |
\lesssim 2^k (2^{n}\delta_1/\delta_0)^{N-j} $ for $1 \le j \le N-1$ and $2^{k-2}/B  \le \big|\big\langle \bar{\mathbf y}_\mu, \mathbf r^{(N)}
  \big\rangle\big|\le CB 2^k$ on $\supp_\xi \bfa_{\nu}^{\mu,n}$.  Thus, as before (cf.  \eqref{supp-supp}),  we  have  
\[ 
\bar{\mathrm Y}_\mu^{\delta_0} ( \supp_{\xi} \bfa_{\nu}^{\mu,n} )
\subset  
\mathbf s\Big(\frac{2^n\delta_1}{\delta_0}\nu-\mu, \, C\frac{2^n\delta_1}{\delta_0},\, CB;  \, \mathbf r_\circ^N \Big) \times \mathbb R^{d-N}
\]
for some $C>0$. Note  $\supp \mathcal F( \mathcal  A[\bfa _{\nu}^{\mu,n}]f)\subset\supp_{\xi}   \bfa_{\nu}^{\mu,n}$. Therefore, changing variables, by \eqref{dgain} with $N$ replaced by $N-1$ and its cylindrical extension (e.g.,\eqref{cylinder}),  we get  
\Be
\label{de-coupn}
\big\| \sum_{\nu\in \mathfrak J_{n}^\mu} \mathcal A[\bfa_{\nu}^{\mu,n}]f \big\|_p
	\le C_0
	\big( {2^n\delta_1} /{\delta_0}\big)^{\frac {N} {p}-1+\epsilon_0}  \big(\sum_{\nu\in \mathfrak J_{n}^\mu} 
	\big\|\mathcal A[\bfa _{\nu}^{\mu,n}]f \big\|_p^p\, \big)^{1/p}
\Ee
 for $2N-2<p <\infty$ (cf. \eqref{decoupn}). Since $\mathcal A [\bfa^\mu]f=\sum_n \sum_{\nu \in \mathfrak J_n^\mu} \mathcal A[\bfa_{\nu}^{\mu,n}]f$,  by Minkowski's inequality and  \eqref{de-coupn},  we have $(\sum_\mu \| \mathcal A[\bfa ^\mu]f\|_p^p\,)^{1/p}$ bounded by 
\begin{equation*}
\begin{aligned} 
 \sum_{n\ge 0} \bar {\mathbf E}_n:=	C_0 \sum_{n\ge 0}
		\big( {2^n\delta_1}/{\delta_0} \big)^{{\frac {N} {p}-1+\epsilon_0}} \big(\sum_{\mu} \sum_{\nu\in \mathfrak J_n^\mu} \| \mathcal A[\bfa ^{\mu,n}_{\nu}]f\|_p^p\,\big)^{1/p}.  
	\end{aligned}
\end{equation*}

The proof of Lemma \ref{itsym}  also shows $C^{-1} \bfa_{\nu}^{\mu,n}\in \bar{\mathfrak A}_k(2^n\delta_1\nu, 2^n\delta_1)$ for a positive constant $C$.  
Therefore, the matter is reduced  to obtaining 
\begin{align}\label{mn1}
 \big(\sum_{\mu} \sum_{\nu \in \mathfrak J_n^\mu} \|\mathcal A[\bfa_{\nu}^{\mu,n}] f\|_{L^p(\mathbb R^d)}^p \,\big)^{1/p}
\lesssim_B\! (2^n\delta_1)^{1-\frac Np} 2^{-\frac kp}
\|f\|_{L^p(\mathbb R^d)}, \qquad n\ge 1
\end{align}
for $p>2(N-2)$. This gives $\sum_{n\ge 1} \bar {\mathbf E}_n \lesssim_B\!  \delta_0^{-N/p+1} 2^{-k/p}\|f\|_p$ since $2^n\delta_1\le C\delta_0$. 

The proof of \eqref{mn1}  is similar with that of  \eqref{T1pp}.    Since $C^{-1} \bfa_{\nu}^{\mu,n}\in  \bar{\mathfrak A}_k(2^n\delta_1\nu, 2^n\delta_1)$,  we  have $\mathcal A[\bfa ^{\mu, n}_{\nu}]f=  \mathcal A[\bfa ^{\mu, n}_{\nu}]   P_{2^n\delta_1\nu}^{2^n\delta_1} f $.  Besides, \eqref{small} or \eqref{lowergg} for some $0\le j\le N-2$ holds on $\supp \bfa_{\nu}^{\mu,n}$. Thus, 
 \eqref{lowD}  holds with $\delta=2^n\delta_1$ for some $B'$ 
on  $\supp  \bfa_{\nu}^{\mu,n}$ for $n\ge 1$ (see \emph{Proof of Lemma \ref{lower-lower}}).   Therefore,  applying Lemma \ref{resc2} to $\mathcal A[\bfa_{\nu}^{\mu,n}] f$ and then the assumption (Theorem \ref{sob-induct} with $L=N-1$), we obtain 
\[ \|\mathcal A[\bfa_{\nu}^{\mu,n}] f\|_{L^p(\mathbb R^d)}\lesssim_B\! (2^n\delta_1)^{1-\frac Np} 2^{-\frac kp}
 \|P_{2^n\delta_1\nu}^{2^n\delta_1} f\|_p.\]   
This combined with \eqref{lplp} gives \eqref{mn1} as desired.

\subsection*{Acknowledgement} H. Ko  was  supported by NRF2019R1A6A3A01092525. S.  Lee and S.  Oh were  supported by  NRF2021R1A2B5B02001786. 
The authors would like to thank  N. Bez, A. Seeger, and J.  Wright for their valuable comments.

\end{document}